\documentclass{article}
\usepackage[utf8]{inputenc}
\usepackage[spanish,english]{babel}
\usepackage{amsmath,amssymb,amsthm,amsfonts}
\usepackage{geometry}
\usepackage{hyperref}
\usepackage{fancyhdr}
\usepackage{titlesec}
\usepackage{listings}
\usepackage{graphicx,graphics}
\usepackage{multicol}
\usepackage{multirow}
\usepackage{color}
\usepackage{float} 
\usepackage{subfig}
\usepackage[figuresright]{rotating}
\usepackage{enumerate}
\usepackage{anysize} 
\usepackage{url}
\title{Notas sobre Teor\'ia de colas y algunas aplicaciones\\
\small{Notes about Queueting Theory and some applications}}
\author{Carlos E. Martínez-Rodríguez\\
	Academia de Matem\'aticas\\
	Universidad Aut\'onoma de la Ciudad de M\'exico\\
	Iztapalapa,  Ermita Iztapalapa 4163\\
	\texttt{carlos.martinez@uacm.edu.mx} }
\date{}

\newtheorem{Def}{Definición}[section]
\newtheorem{Ejem}{Ejemplo}[section]
\newtheorem{Teo}{Teorema}[section]
\newtheorem{Note}{Nota}[section]
\newtheorem{Prop}{Proposición}[section]
\newtheorem{Cor}{Corolario}[section]
\newtheorem{Lema}{Lema}[section]

\newcommand{\nat}{\mathbb{N}}
\newcommand{\ent}{\mathbb{Z}}
\newcommand{\rea}{\mathbb{R}}
\newcommand{\Eb}{\mathbf{E}}
\newcommand{\esp}{\mathbb{E}}
\newcommand{\prob}{\mathbb{P}}
\newcommand{\indora}{\mbox{$1$\hspace{-0.8ex}$1$}}

\renewcommand{\abstractname}{Resumen}
\numberwithin{equation}{section}

\begin{document}
\maketitle

\begin{abstract}
En este trabajo se presenta una revisión exhaustiva de los procesos estocásticos, con un enfoque particular en las cadenas de Markov y los procesos de salto. Se analizan los principales resultados relacionados con los sistemas de espera. Además, se exponen condiciones que garantizan la estabilidad, o ergodicidad, de dichos sistemas. Asimismo, se discuten los resultados de estabilidad en redes de colas y su extensión a los sistemas de visitas. Finalmente, se presentan contribuciones clave relacionadas con la Función Generadora de Probabilidades, una herramienta esencial en el análisis de los procesos previamente mencionados. La revisión se lleva a cabo desde la perspectiva de la teoría de colas, fundamentada en la notación de Kendall-Lee, subrayando los resultados sobre estabilidad y el cálculo de medidas de desempeño en función de las características particulares de cada proceso.

\end{abstract}

\begin{otherlanguage}{english}
\renewcommand{\abstractname}{Abstract} 
\begin{abstract}
This paper presents a comprehensive review of stochastic processes, with a particular focus on Markov chains and jump processes. The main results related to queuing systems are analyzed. Additionally, conditions that ensure the stability, or ergodicity, of such systems are presented. The paper also discusses stability results for queuing networks and their extension to visiting systems. Finally, key contributions concerning the Probability Generating Function, an essential tool in the analysis of the aforementioned processes, are introduced. The review is conducted from the perspective of queuing theory, grounded in the Kendall-Lee notation, emphasizing stability results and the computation of performance measures based on the specific characteristics of each process.
\end{abstract}
\end{otherlanguage}


\section*{Introducción}
Los sistemas de espera consisten en una o m\'as colas en las que los usuarios llegan para ser atendidos conforme a una pol\'itica de servicio espec\'ifica \cite{Semenova} . Los arribos siguen una ley de distribuci\'on por unidad o intervalo de tiempo, lo que implica que uno de los procesos clave a modelar son los tiempos entre llegadas y, por lo tanto, la distribuci\'on asociada a estos. En la mayor\'ia de los casos, las estaciones de servicio o colas operan bajo la pol\'itica \textit{first in, first out} (FIFO) \cite{Kleinrock}, donde el primer usuario en llegar es el primero en ser atendido. Otra caracter\'istica relevante de los sistemas de espera es la capacidad de las colas, la cual puede ser finita o infinita, permitiendo que un n\'umero ilimitado de usuarios ingrese a la cola y sea atendido hasta vaciarla por completo \cite{Asmussen}, \cite{CooperI}, \cite{LevySidi}. Una vez que el sistema ha sido caracterizado, resulta esencial determinar sus principales m\'etricas de desempe\~no, tales como los tiempos promedio de servicio, atenci\'on y espera. Estas tres medidas son fundamentales para caracterizar el sistema y posibilitan el c\'alculo de las medidas de desempe\~no estacionarias \cite{Roubos, TakagiI}.

El uso de una notaci\'on adecuada que facilite la comprensi\'on de los procesos involucrados en el sistema de espera es crucial, ya que permite identificar si los tiempos de llegada o de servicio son generales, exponenciales, uniformes o deterministas ($G, M, U, D$), o si la capacidad de la cola es finita ($K$) o infinita ($\infty$). Tambi\'en define si el sistema est\'a compuesto por uno o varios servidores, y finalmente, especifica la disciplina de servicio que siguen los servidores para atender la cola \cite{TakagiI}.

Por otro lado, las redes de colas \cite{BosBoon} constituyen una generalizaci\'on de los sistemas de colas, donde es posible la existencia de m\'ultiples colas atendidas de manera independiente por uno o varios servidores. En este contexto, cada cola posee sus propios servidores que atienden a los usuarios hasta que la cola se vac\'ia por completo. Las redes de colas tienen aplicaciones naturales en diversos contextos, como cajas de pago en supermercados, centros de atenci\'on telef\'onica (\textit{call centers}), procesos de manufactura, entre otros \cite{CooperI, Bhat}.

Los sistemas de visitas \cite{Semenova} representan una extensi\'on de las redes de colas, cuya principal diferencia radica en que las colas se encuentran ubicadas en proximidad, y uno o varios servidores se desplazan entre ellas para proporcionar servicio. Estos desplazamientos pueden ser c\'iclicos, peri\'odicos, deterministas o aleatorios. Los servidores pueden seguir diferentes pol\'iticas de servicio, como atender \'unicamente a los usuarios presentes al momento de su llegada o tambi\'en a aquellos que llegan mientras se est\'a proporcionando el servicio. Una vez que la cola se vac\'ia, los servidores se trasladan a la siguiente cola y comienzan a atender conforme a la pol\'itica establecida en ese sistema. Un aspecto relevante de los sistemas de visitas es el tiempo de traslado entre colas, lo cual permite obtener m\'etricas de desempe\~no como la longitud promedio de la cola o el tiempo promedio de servicio en cada estaci\'on.

La revisi\'on de estos temas se llev\'o a cabo en su momento bajo la minuciosa y puntual supervisi\'on del Dr. Ra\'ul Montes de Oca Machorro y la Dra. Patricia Saavedra Barrera, cuyas oportunas y valiosas sugerencias y comentarios fueron fundamentales para desarrollar el estudio de este tipo de procesos estoc\'asticos.

Este documento no pretende ser un estudio exhaustivo de los sistemas de visitas, sino m\'as bien un an\'alisis que recupere los elementos fundamentales para estudiar y comprender este tipo particular de procesos estoc\'asticos. El objetivo de este trabajo es facilitar el acercamiento y estudio de procesos cuya aplicaci\'on es directa y de gran utilidad en diversas \'areas \cite{Semenova,LevySidi}.

El documento est\'a organizado de la siguiente manera: en primer lugar, se presenta una revisi\'on general de la teor\'ia subyacente a los procesos estoc\'asticos. Posteriormente, se exploran procesos particularmente relevantes como los procesos de Markov de saltos, destacando los resultados b\'asicos que, a juicio del autor, son fundamentales para introducirse en el tema. Una vez establecidos los procesos de Markov de saltos, se introduce la notaci\'on de Kendall-Lee para facilitar el estudio de las colas, presentando una amplia variedad de resultados sobre las principales medidas de desempe\~no y condiciones de ergodicidad. En la Secci\'on 4 se realiza una breve menci\'on a las redes de colas abiertas. Finalmente, en la Secci\'on 5 se describen los sistemas de visitas, destacando los elementos principales que los definen y los resultados m\'as relevantes sobre tiempos de espera promedio y longitud promedio de las colas, basados en la obra de Takagi. Adicionalmente, se presentan resultados \'utiles sobre la funci\'on generadora de probabilidades, una herramienta clave en el an\'alisis de los sistemas de visitas.

\section{Procesos Estoc\'asticos}\label{Procesos.Estocasticos}

\begin{Def}\index{Conjunto Medible}
Sea $X$ un conjunto y $\mathcal{F}$ una $\sigma$-\'algebra de subconjuntos de $X$, la pareja $\left(X,\mathcal{F}\right)$ es llamado espacio medible. Un subconjunto $A$ de $X$ es llamado medible, o medible con respecto a $\mathcal{F}$, si $A\in\mathcal{F}$.
\end{Def}

\begin{Def}\index{Medida $\sigma$-finita}
Sea $\left(X,\mathcal{F},\mu\right)$ espacio de medida. Se dice que la medida $\mu$ es $\sigma$-finita si se puede escribir $X=\bigcup_{n\geq1}X_{n}$ con $X_{n}\in\mathcal{F}$ y $\mu\left(X_{n}\right)<\infty$.
\end{Def}

\begin{Def}\label{Cto.Borel}\index{Conjunto de Borel}
Sea $X$ un espacio topológico. El álgebra de Borel en $X$, denotada por $\mathcal{B}(X)$, es la $\sigma$-álgebra generada por la colección de todos los conjuntos abiertos de $X$. Es decir, $\mathcal{B}(X)$ es la colección más pequeña de subconjuntos de $X$ que contiene todos los conjuntos abiertos y es cerrada bajo la unión numerable, la intersección numerable y el complemento.
\end{Def}

\begin{Def}\label{Funcion.Medible}\index{Funci\'on Medible}
Una funci\'on $f:X\rightarrow\rea$, es medible si para cualquier n\'umero real $\alpha$ el conjunto \[\left\{x\in X:f\left(x\right)>\alpha\right\},\] pertenece a $X$. Equivalentemente, se dice que $f$ es medible si \[f^{-1}\left(\left(\alpha,\infty\right)\right)=\left\{x\in X:f\left(x\right)>\alpha\right\}\in\mathcal{F}.\]
\end{Def}

\begin{Def}\label{Def.Cilindros}\index{Cilindro}
Sean $\left(\Omega_{i},\mathcal{F}_{i}\right)$, $i=1,2,\ldots,$ espacios medibles y $\Omega=\prod_{i=1}^{\infty}\Omega_{i}$ el conjunto de todas las sucesiones $\left(\omega_{1},\omega_{2},\ldots,\right)$ tales que $\omega_{i}\in\Omega_{i}$, $i=1,2,\ldots,$. Si $B^{n}\subset\prod_{i=1}^{\infty}\Omega_{i}$, definimos $B_{n}=\left\{\omega\in\Omega:\left(\omega_{1},\omega_{2},\ldots,\omega_{n}\right)\in B^{n}\right\}$. Al conjunto $B_{n}$ se le llama {\em cilindro} con base $B^{n}$, el cilindro es llamado medible si $B^{n}\in\prod_{i=1}^{\infty}\mathcal{F}_{i}$.
\end{Def}

\begin{Def}\label{Def.Proc.Adaptado}\index{Proceso Adaptado}
Sea $X\left(t\right),t\geq0$ proceso estoc\'astico, el proceso es adaptado a la familia de $\sigma$-\'algebras $\mathcal{F}_{t}$, para $t\geq0$, si para $s<t$ implica que $\mathcal{F}_{s}\subset\mathcal{F}_{t}$, y $X\left(t\right)$ es $\mathcal{F}_{t}$-medible para cada $t$. Si no se especifica $\mathcal{F}_{t}$ entonces se toma $\mathcal{F}_{t}$ como $\mathcal{F}\left(X\left(s\right),s\leq t\right)$, la m\'as peque\~na $\sigma$-\'algebra de subconjuntos de $\Omega$ que hace que cada $X\left(s\right)$, con $s\leq t$ sea Borel medible.
\end{Def}

\begin{Def}\label{Def.Tiempo.Paro}\index{Tiempos de Paro}
Sea $\left\{\mathcal{F}\left(t\right),t\geq0\right\}$ familia creciente de sub $\sigma$-\'algebras. es decir, $\mathcal{F}\left(s\right)\subset\mathcal{F}\left(t\right)$ para $s\leq t$. Un tiempo de paro para $\mathcal{F}\left(t\right)$ es una funci\'on $T:\Omega\rightarrow\left[0,\infty\right]$ tal que $\left\{T\leq t\right\}\in\mathcal{F}\left(t\right)$ para cada $t\geq0$. Un tiempo de paro para el proceso estoc\'astico $X\left(t\right),t\geq0$ es un tiempo de paro para las $\sigma$-\'algebras $\mathcal{F}\left(t\right)=\mathcal{F}\left(X\left(s\right)\right)$.
\end{Def}

\begin{Def}\index{Proceso Adaptado}
Sea $X\left(t\right),t\geq0$ proceso estoc\'astico, con $\left(S,\chi\right)$ espacio de estados. Se dice que el proceso es adaptado a $\left\{\mathcal{F}\left(t\right)\right\}$, es decir, si para cualquier $s,t\in I$, $I$ conjunto de \'indices, $s<t$, se tiene que $\mathcal{F}\left(s\right)\subset\mathcal{F}\left(t\right)$, y $X\left(t\right)$ es $\mathcal{F}\left(t\right)$-medible,
\end{Def}

\begin{Def}\index{Proceso de Markov}
Sea $X\left(t\right),t\geq0$ proceso estoc\'astico, se dice que es un Proceso de Markov relativo a $\mathcal{F}\left(t\right)$ o que $\left\{X\left(t\right),\mathcal{F}\left(t\right)\right\}$ es de Markov si y s\'olo si para cualquier conjunto $B\in\chi$,  y $s,t\in I$, $s<t$ se cumple que
\begin{equation}\label{Prop.Markov}
P\left\{X\left(t\right)\in B|\mathcal{F}\left(s\right)\right\}=P\left\{X\left(t\right)\in B|X\left(s\right)\right\}.
\end{equation}
\end{Def}

\begin{Note}
Si se dice que $\left\{X\left(t\right)\right\}$ es un Proceso de Markov sin mencionar $\mathcal{F}\left(t\right)$, se asumir\'a que 
\begin{eqnarray*}
\mathcal{F}\left(t\right)=\mathcal{F}_{0}\left(t\right)=\mathcal{F}\left(X\left(r\right),r\leq t\right),
\end{eqnarray*}
entonces la ecuaci\'on (\ref{Prop.Markov}) se puede escribir como
\begin{equation}
P\left\{X\left(t\right)\in B|X\left(r\right),r\leq s\right\} = P\left\{X\left(t\right)\in B|X\left(s\right)\right\}.
\end{equation}
\end{Note}

\begin{Teo}
Sea $\left(X_{n},\mathcal{F}_{n},n=0,1,\ldots,\right\}$ Proceso de Markov con espacio de estados $\left(S_{0},\chi_{0}\right)$ generado por una distribuici\'on inicial $P_{o}$ y probabilidad de transici\'on $p_{mn}$, para $m,n=0,1,\ldots,$ $m<n$, que por notaci\'on se escribir\'a como $p\left(m,n,x,B\right)\rightarrow p_{mn}\left(x,B\right)$. Sea $S$ tiempo de paro relativo a la $\sigma$-\'algebra $\mathcal{F}_{n}$. Sea $T$ funci\'on medible, $T:\Omega\rightarrow\left\{0,1,\ldots,\right\}$. Sup\'ongase que $T\geq S$, entonces $T$ es tiempo de paro. Si $B\in\chi_{0}$,
entonces
\begin{equation}\label{Prop.Fuerte.Markov}
P\left\{X\left(T\right)\in B,T<\infty|\mathcal{F}\left(S\right)\right\} = p\left(S,T,X\left(s\right),B\right).
\end{equation}
en $\left\{T<\infty\right\}$.
\end{Teo}

\subsection*{Cadenas de Markov}

\begin{Def}\index{Cadena de Markov}
Sea $\left(\Omega,\mathcal{F},\prob\right)$ un espacio de probabilidad y $\mathbf{E}$ un conjunto no vac\'io, finito o numerable. Una sucesi\'on de variables aleatorias $\left\{X_{n}:\Omega\rightarrow\mathbf{E},n\geq0\right\}$ se le llama \textit{Cadena de Markov} con espacio de estados $\mathbf{E}$ si satisface la condici\'on de Markov, esto es, si para todo $n\geq1$ y toda sucesi\'on $x_{0},x_{1},\ldots,x_{n},x,y\in\mathbf{E}$ se cumple que 

\begin{equation}
P\left\{X_{n}=y|X_{n-1}=x,\ldots,X_{0}=x_{0}\right\}=P\left\{X_{n}=x_{n}|X_{n-1}=x_{n-1}\right\}.
\end{equation}
La distribuci\'on de $X_{0}$ se llama distribuci\'on inicial y se denotar\'a por $\pi$.
\end{Def}

\begin{Note}\index{Probabilidades Condicionales}
Las probabilidades condicionales $P\left\{X_{n}=y|X_{n-1}=x\right\}$ se les llama \textit{probabilidades condicionales}\index{Probabilidades Condicionales}
\end{Note}

\begin{Note}\index{Cadenas Homog\'eneas}
En este trabajo se considerar\'an solamente aquellas cadenas de Markov con probabilidades de transici\'on estacionarias, es decir, aquellas que no dependen del valor de $n$ (se dice que es una cadena homog\'enea), es decir, cuando se diga $X_{n},n\geq0$ es cadena de Markov, se entiende que es una sucesi\'on de variables aleatorias que satisfacen la propiedad de Markov y que tienen probabilidades de transici\'on estacionarias.\index{Cadena Homog\'enea}
\end{Note}

\begin{Note}
Para una cadena de Markov Homog\'enea se tiene la siguiente denotaci\'on
\begin{equation}
P\left\{X_{n}=y|X_{n-1}=x\right\}=P_{x,y}.
\end{equation}
\end{Note}

\begin{Note}\index{Probabilidades de Transici\'on}
Para $m\geq1$ se denotar\'a por $P^{(m)}_{x,y}$ a $P\left\{X_{n+m}=y|X_{n}=x\right\}$, que significa la probabilidad de ir en $m$ pasos o unidades de tiempo de $x$ a $y$, y se le llama \textit{probabilidad de transici\'on en $m$ pasos}.
\end{Note}

\begin{Note}\index{Delta de Kronecker}
Para $x,y\in\mathbf{E}$ se define a $P^{(0)}_{x,y}$ como $\delta_{x,y}$, donde $\delta_{x,y}$ es la delta de Kronecker, es decir, vale 1 si $x=y$ y 0 en otro caso.
\end{Note}

\begin{Note}\index{Matriz de Transici\'on}
En el caso de que $\mathbf{E}$ sea finito, se considera la matrix $P=\left(P_{x,y}\right)_{x,y\in \mathbf{E}}$ y se le llama \textit{matriz de transici\'on}.
\end{Note}

\begin{Note}
Si la distribuci\'on inicial $\pi$ es igual al vector $\left(\delta_{x,y}\right)_{y\in\mathbf{E}}$, es decir,
\begin{eqnarray*}
P\left(X_{0}=x)=1\right)\textrm{ y }P\left(X_{0}\neq x\right)=0,
\end{eqnarray*}
entonces se toma la notaci\'on 
\begin{eqnarray}
&&P_{x}\left(A\right)=P\left(A|X_{0}=x\right),A\in\mathcal{F},
\end{eqnarray}
y se dice que la cadena empieza en $A$. Se puede demostrar que $P_{x}$ es una nueva medida de probabilidad en el espacio $\left(\Omega,\mathcal{F}\right)$.
\end{Note}

\begin{Note}
La suma de las entradas de los renglones de la matriz de transici\'on es igual a uno, es decir, para todo $x\in \mathbf{E}$ se tiene $\sum_{y\in\mathbf{E}}P_{x,y}=1$.
\end{Note}

Para poder obtener uno de los resultados m\'as importantes en cadenas de Markov, la \textit{ecuaci\'on de Chapman-kolmogorov} se requieren los siguientes resultados:

\begin{Lema}
Sean $x,y,z\in\Eb$ y $0\leq m\leq n-1$, entonces se cumple que
\begin{equation}
P\left(X_{n+1}=y|X_{n}=z,X_{m}=x\right)=P_{z,y}.
\end{equation}
\end{Lema}

\begin{Prop}
Si $x_{0},x_{1},\ldots,x_{n}\in \Eb$ y $\pi\left(x_{0}\right)=P\left(X_{0}=x_{0}\right)$, entonces
\begin{equation}
P\left(X_{1}=x_{1},\ldots,X_{n}=x_{n},X_{0}=x_{0}\right)=\pi\left(x_{0}\right)P_{x_{0},x_{1}}\cdot P_{x_{1},x_{2}}\cdots P_{x_{n-1},x_{n}}.
\end{equation}
\end{Prop}

De la proposici\'on anterior se tiene
\begin{equation}
P\left(X_{1}=x_{1},\ldots,X_{n}=x_{n}|X_{0}=x_{0}\right)=P_{x_{0},x_{1}}\cdot P_{x_{1},x_{2}}\cdots P_{x_{n-1},x_{n}}.
\end{equation}

finalmente tenemos la siguiente proposici\'on:

\begin{Prop}
Sean $n,k\in\nat$ fijos y $x_{0},x_{1},\ldots,x_{n},\ldots,x_{n+k}\in\Eb$, entonces
\begin{eqnarray*}
&&P\left(X_{n+1}=x_{n+1},\ldots,X_{n+k}=x_{n+k}|X_{n}=x_{n},\ldots,X_{0}=x_{0}\right)\\
&=&P\left(X_{1}=x_{n+1},X_{2}=x_{n+2},\cdots,X_{k}=x_{n+k}|X_{0}=x_{n}\right).
\end{eqnarray*}
\end{Prop}

\begin{Ejem}
Sea $X_{n}$ una variable aleatoria al tiempo $n$ tal que
\begin{eqnarray}
\begin{array}{l}
P\left(X_{n+1}=1 \mid X_{n}=0\right)=p,\\
P\left(X_{n+1}=0 \mid X_{n}=1\right)=q=1-p,\\
P\left(X_{0}=0\right)=\pi_{0}\left(0\right).
\end{array}
\end{eqnarray}

\end{Ejem}

Se puede demostrar que
\begin{eqnarray}
\begin{array}{l}
P\left(X_{n}=0\right)=\frac{q}{p+q},\\
P\left(X_{n}=1\right)=\frac{p}{p+q}.
\end{array}
\end{eqnarray}

\begin{Ejem}
El problema de la Caminata Aleatoria.
\end{Ejem}

\begin{Ejem}
El problema de la ruina del jugador.
\end{Ejem}

\begin{Ejem}
Sea $\left\{Y_{i}\right\}_{i=0}^{\infty}$ sucesi\'on de variables aleatorias independientes e identicamente distribuidas, definidas sobre un espacio de probabilidad $\left(\Omega,\mathcal{F},\prob\right)$ y que toman valores enteros, se tiene que la sucesi\'on $\left\{X_{i}\right\}_{i=0}^{\infty}$ definida por $X_{j}=\sum_{i=0}^{j}Y_{i}$ es una cadena de Markov en el conjunto de los n\'umeros enteros.
\end{Ejem}

\begin{Prop}\index{Ecuaciones de Chapman-Kolmogorov}
Para una cadena de Markov $\left(X_{n}\right)_{n\in\nat}$ con espacio de estados $\Eb$ y para todo $n,m\in \nat$ y toda pareja $x,y\in\Eb$ se cumple
\begin{equation}
P\left(X_{n+m}=y|X_{0}=x\right)=\sum_{z\in\Eb}P_{x,z}^{(m)}P_{z,y}^{(n)}=P_{x,y}^{(n+m)}.
\end{equation}
\end{Prop}

\begin{Note}
Para una cadena de Markov con un n\'umero finito de estados, se puede pensar a $P^{n}$ como la $n$-\'esima potencia de la matriz $P$. Sea $\pi_{0}$ distribuci\'on inicial de la cadena de Markov, como 
\begin{eqnarray}
P\left(X_{n}=y\right)=\sum_{x} P\left(X_{0}=x,X_{n}=y\right)=\sum_{x} P\left(X_{0}=x\right)P\left(X_{n}=y|X_{0}=x\right),
\end{eqnarray}
se puede comprobar que 

\begin{eqnarray}
P\left(X_{n}=y\right)=\sum_{x} \pi_{0}\left(X\right)P^{n}\left(x,y\right).
\end{eqnarray}
\end{Note}

Con lo anterior es posible calcular la distribuici\'on de $X_{n}$ en t\'erminos de la distribuci\'on inicial $\pi_{0}$ y la funci\'on de transici\'on de $n$-pasos $P^{n}$,
\begin{eqnarray}
P\left(X_{n+1}=y\right)=\sum_{x} P\left(X_{n}=x\right)P\left(x,y\right).
\end{eqnarray}
\begin{Note}
Si se conoce la distribuci\'on de $X_{0}$ se puede conocer la distribuci\'on de $X_{1}$.
\end{Note}

\subsection*{Clasificaci\'on de Estados}

\begin{Def}\index{Tiempos de Paro}
Para $A$ conjunto en el espacio de estados, se define un tiempo de paro $T_{A}$ de $A$ como
\begin{equation}
T_{A}=min_{n>0}\left(X_{n}\in A\right).
\end{equation}
\end{Def}

\begin{Note}
Si $X_{n}\notin A$ para toda $n>0$, $T_{A}=\infty$, es decir,  $T_{A}$ es el primer tiempo positivo que la cadena de Markov est\'a en $A$.
\end{Note}

Una vez que se tiene la definici\'on anterior se puede demostrar la siguiente igualdad:

\begin{Prop}
$P^{n}\left(x,y\right)=\sum_{m=1}^{n}P_{x}\left(T_{y}=m\right)P^{n.m}\left(y,x\right), n\geq1$.
\end{Prop}
\medskip

\begin{Def}
En una cadena de Markov $\left(X_{n}\right)_{n\in\nat}$ con espacio de estados $\Eb$, matriz de transici\'on $\left(P_{x,y}\right)_{x,y\in\Eb}$ y para $x,y\in\Eb$,  se dice que
\begin{itemize}
\item[a) ]  De $x$ se accede a $y$ si existe $n\geq0$ tal que $P_{x,y}^{(n)}>0$ y se denota por $\left(x\rightarrow y\right)$.

\item[b) ] $x$ y $y$ se comunican entre s\'i, lo que se denota por $\left(x\leftrightarrow y\right)$, si se cumplen $\left(x\rightarrow y\right)$ y $\left(y\rightarrow x\right)$.

\item[c) ] Un estado $x\in\Eb$ es estado recurrente si $$P\left(X_{n}=x\textrm{ para alg\'un }n\in\nat|X_{0}=x \right)\equiv1.$$ \index{Estados recurrentes}

\item[d) ] Un estado $x\in\Eb$ es estado transitorio si $$P\left(X_{n}=x\textrm{ para alg\'un }n\in\nat|X_{0}=x \right)<1.$$ \index{Estados transitorios}

\item[e) ] Un estado $x\in\Eb$ se llama absorbente si $P_{x,x}\equiv1$.\index{Estados absorbentes}
\end{itemize}
\end{Def}

Se tiene el siguiente resultado:

\begin{Prop}
$x\leftrightarrow y$ es una relaci\'on de equivalencia y da lugar a una partici\'on del espacio de estados $\Eb$.
\end{Prop}

\begin{Def}
Para $E$ espacio de estados
\begin{itemize}

\item[a)  ] Se dice que $C\subset \Eb$ es una clase de comunicaci\'on si cualesquiera dos estados de $C$ se comunic\'an entre s\'i.\index{Clases de Comunicaci\'on}

\item[b)  ] Dado $x\in\Eb$, su clase de comunicaci\'on se denota por: $C\left(x\right)=\left\{y\in\Eb:x\leftrightarrow y\right\}$.

\item[c)  ] Se dice que un conjunto de estados  $C\subset \Eb$ es cerrado si ning\'un estado de $\Eb-C$ puede ser accedido desde un estado de $C$.
\end{itemize}
\end{Def}

\begin{Def}\index{Cadena Irreducible}
Sea $\Eb$ espacio de estados, se dice que la cadena es irreducible si cualquiera de las siguientes condiciones, equivalentes entre s\'i,  se cumplen
\begin{enumerate}
\item[a) ] Desde cualquier estado de $\Eb$ se puede acceder a cualquier otro.

\item[b) ] Todos los estados se comunican entre s\'i.

\item[c) ] $C\left(x\right)=\Eb$ para alg\'un $x\in\Eb$.

\item[d) ] $C\left(x\right)=\Eb$ para todo $x\in\Eb$.

\item[e) ] El \'unico conjunto cerrado es el total.
\end{enumerate}
\end{Def}
Por lo tanto tenemos la siguiente proposici\'on:
\begin{Prop}  Sea $\Eb$ espacio de estados y $T$ tiempo de paro, entonces se tiene que
\begin{enumerate}
\item[a) ] Un estado $x\in\Eb$ es recurrente si y s\'olo si $P\left(T_{x}<\infty|x_{0}=x\right)=1$.

\item[b) ] Un estado $x\in\Eb$ es transitorio si y s\'olo si $P\left(T_{x}<\infty|x_{0}=x\right)<1$.

\item[c) ] Un estado $x\in\Eb$ es absorbente si y s\'olo si $P\left(T_{x}=1|x_{0}=x\right)=1$.

\end{enumerate}
\end{Prop}

Sea $v=\left(v_{i}\right)_{i\in E}$ medida no negativa en $E$, podemos definir una nueva medida $v\prob$ que asigna masa $\sum_{i\in E}v_{i}p_{ij}$ a cada estado $j$.

\begin{Def}\index{Medida Estacionaria}
La medida $v$ es estacionaria si $v_{i}<\infty$ para toda $i$ y adem\'as $v\prob=v$.
\end{Def}
En el caso de que $v$ sea distribuci\'on, independientemente de que sea estacionaria o no, se cumple con

\begin{eqnarray}
\prob_{v}\left[X_{1}=j\right]=\sum_{i\in E}\prob_{v}\left[X_{0}=i\right]p_{ij}=\sum_{i\in E}v_{i}p_{ij}=\left(vP\right)_{j}.
\end{eqnarray}

\begin{Teo}
Supongamos que $v$ es una distribuci\'on estacionaria. Entonces
\begin{itemize}
\item[i)] La cadena es estrictamente estacionaria con respecto a
$\prob_{v}$, es decir, $\prob_{v}$-distribuci\'on de $\left\{X_{n},X_{n+1},\ldots\right\}$ no depende de $n$;
\item[ii)] Existe un aversi\'on estrictamente estacionaria $\left\{X_{n}\right\}_{n\in Z}$ de la cadena con doble tiempo infinito y $\prob\left(X_{n}=i\right)=v_{i}$ para toda $n\in Z$.
\end{itemize}
\end{Teo}

\begin{Teo}
Sea $i$ estado fijo, recurrente. Entonces una medida estacionaria $v$ puede definirse haciendo que $v_{j}$ sea el n\'umero esperado de visitas a $j$ entre dos visitas consecutivas $i$,

\begin{equation}\label{Eq.3.1}
v_{j}=\esp_{i}\sum_{n=0}^{\tau(i)-1}\indora\left(X_{n}=i\right)=\sum_{n=0}^{\infty}\prob_{i}\left[X_{n}=j,\tau(i)>n\right].
\end{equation}
\end{Teo}

\begin{Teo}\label{Teo.3.3}
Si la cadena es irreducible y recurrente, entonces existe una medida estacionaria $v$, tal que satisface $0<v_{j}<\infty$ para toda $j$, y es \'unica salvo factores multiplicativos, es decir, si $v,v^{*}$ son estacionarias, entonces $c=cv^{*}$ para alguna $c\in\left(0,\infty\right)$.
\end{Teo}

\begin{Cor}\label{Cor.3.5}
Si la cadena es irreducible y positiva recurrente, existe una \'unica distribuci\'on estacionaria $\pi$ dada por
\begin{equation}
\pi_{j}=\frac{1}{\esp_{i}\tau_{i}}\esp_{i}\sum_{n=0}^{\tau\left(i\right)-1}\indora\left(X_{n}=j\right)=\frac{1}{\esp_{j}\tau\left(j\right)}.
\end{equation}
\end{Cor}

\begin{Cor}\label{Cor.3.6}\index{Cadena Positiva Recurrente}
Cualquier cadena de Markov irreducible con un espacio de estados finito es positiva recurrente.
\end{Cor}

\begin{Def}\label{Def.Armonica}\index{Funci\'on Arm\'onica}
Una funci\'on Arm\'onica es el eigenvector derecho $h$ de $P$ correspondiente al eigenvalor 1.
\end{Def}
\begin{eqnarray}
Ph=h\Leftrightarrow h\left(i\right)=\sum_{j\in E}p_{ij}h\left(j\right)=\esp_{i}h\left(X_{1}\right)=\esp\left[h\left(X_{n+1}\right)|X_{n}=i\right].
\end{eqnarray}
es decir, $\left\{h\left(X_{n}\right)\right\}$ es martingala.\\

\begin{Prop}\label{Prop.5.2}\index{Cadena Transitoria}
Sea $\left\{X_{n}\right\}$ cadena irreducible  y sea $i$ estado fijo arbitrario. Entonces la cadena es transitoria s\'i y s\'olo si existe una funci\'on no cero, acotada $h:E-\left\{i\right\}\rightarrow\rea$ que satisface
\begin{equation}\label{Eq.5.1}
h\left(j\right)=\sum_{k\neq i}p_{jk}h\left(k\right)\textrm{   para }j\neq i.
\end{equation}
\end{Prop}

\begin{Prop}\label{Prop.5.4}
Suponga que la cadena es irreducible y sea $E_{0}$ un subconjunto finito de $E$ tal que se cumple la ecuaci\'on \ref{Eq.5.1} para alguna funci\'on $h$ acotada que satisface $h\left(i\right)<h\left(j\right)$ para alg\'un estado $i\notin E_{0}$ y todo $j\in E_{0}$. Entonces la cadena es transitoria.
\end{Prop}

\begin{Lema}\index{Cadena Positiva Recurrente}
Sea $\left\{X_{n}\right\}$ cadena irreducible y se $F$ subconjunto finito del espacio de estados. Entonces la cadena es positiva recurrente si $\esp_{i}\tau\left(F\right)<\infty$ para todo $i\in F$.
\end{Lema}

\begin{Prop}
Sea $\left\{X_{n}\right\}$ cadena irreducible y transiente o cero recurrente, entonces $p_{ij}^{n}\rightarrow0$ conforme $n\rightarrow\infty$ para cualquier $i,j\in E$, $E$ espacio de estados.
\end{Prop}

Se tiene el siguiente resultado:

\begin{Teo}
Sea $\left\{X_{n}\right\}$ cadena irreducible y aperi\'odica positiva recurrente, y sea $\pi=\left\{\pi_{j}\right\}_{j\in E}$ la distribuci\'on estacionaria. Entonces $p_{ij}^{n}\rightarrow\pi_{j}$ para todo $i,j$.
\end{Teo}

\begin{Def}\label{Def.Ergodicidad}\index{Cadena Erg\'odica}
Una cadena irreducible aperiodica, positiva recurrente con medida estacionaria $v$, es llamada {\em erg\'odica}.
\end{Def}

\begin{Prop}\label{Prop.4.4}
Sea $\left\{X_{n}\right\}$ cadena irreducible y recurrente con medida estacionaria $v$, entonces para todo $i,j,k,l\in E$
\begin{equation}
\frac{\sum_{n=0}^{m}p_{ij}^{n}}{\sum_{n=0}^{m}p_{lk}^{n}}\rightarrow\frac{v_{j}}{v_{k}}\textrm{,    }m\rightarrow\infty
\end{equation}
\end{Prop}

\begin{Lema}\label{Lema.4.5}
La matriz $\widetilde{P}$ con elementos 
\begin{eqnarray}
\widetilde{p}_{ij}=\frac{v_{ji}p_{ji}}{v_{i}}
\end{eqnarray}
es una matriz de transici\'on. Adem\'s, el $i$-\'esimo elementos $\widetilde{p}_{ij}^{m}$ de la matriz potencia $\widetilde{P}^{m}$ est\'a dada por 

\begin{eqnarray}
\widetilde{p}_{ij}^{m}=\frac{v_{ji}p_{ji}^{m}}{v_{i}}.
\end{eqnarray}
\end{Lema}

\begin{Lema}
Def\'inase 
\begin{eqnarray}
N_{i}^{m}=\sum_{n=0}^{m}\indora\left(X_{n}=i\right)
\end{eqnarray} 
como el n\'umero de visitas a $i$ antes del tiempo $m$. Entonces si la cadena es reducible y recurrente, 
\begin{eqnarray}
lim_{m\rightarrow\infty}\frac{\esp_{j}N_{i}^{m}}{\esp_{k}N_{i}^{m}}=1\textrm{ para todo }j,k\in E.
\end{eqnarray}
\end{Lema}

\subsection*{Ejemplos}

Supongamos que se tiene la siguiente cadena:
\begin{equation}
\left(\begin{array}{cc}
1-q & q\\
p & 1-p\\
\end{array}
\right).
\end{equation}
Si $P\left[X_{0}=0\right]=\pi_{0}(0)=a$ y $P\left[X_{0}=1\right]=\pi_{0}(1)=b=1-\pi_{0}(0)$, con $a+b=1$, entonces despu\'es de un procedimiento m\'as o menos corto se tiene que:

\begin{eqnarray*}
P\left[X_{n}=0\right]=\frac{p}{p+q}+\left(1-p-q\right)^{n}\left(a-\frac{p}{p+q}\right).\\
P\left[X_{n}=1\right]=\frac{q}{p+q}+\left(1-p-q\right)^{n}\left(b-\frac{q}{p+q}\right).\\
\end{eqnarray*}
donde, como $0<p,q<1$, se tiene que $|1-p-q|<1$, entonces $\left(1-p-q\right)^{n}\rightarrow 0$ cuando $n\rightarrow\infty$. Por lo tanto
\begin{eqnarray*}
lim_{n\rightarrow\infty}P\left[X_{n}=0\right]=\frac{p}{p+q}.\\
lim_{n\rightarrow\infty}P\left[X_{n}=1\right]=\frac{q}{p+q}.
\end{eqnarray*}
Si hacemos $v=\left(\frac{p}{p+q},\frac{q}{p+q}\right)$, entonces
\begin{eqnarray*}
\left(\frac{p}{p+q},\frac{q}{p+q}\right)\left(\begin{array}{cc}
1-q & q\\
p & 1-p\\
\end{array}\right).
\end{eqnarray*}

\begin{Prop}\label{Prop.5.4}
Suponga que la cadena es irreducible y sea $E_{0}$ un subconjunto finito de $E$ tal que se cumple la ecuaci\'on \ref{Eq.5.1} para alguna funci\'on $h$ acotada que satisface $h\left(i\right)<h\left(j\right)$ para alg\'un estado $i\notin E_{0}$ y todo $j\in E_{0}$. Entonces la cadena es transitoria.
\end{Prop}

\section{Procesos de Markov de Saltos}
Consideremos un estado que comienza en el estado $x_{0}$ al tiempo $0$, supongamos que el sistema permanece en $x_{0}$ hasta alg\'un tiempo positivo $\tau_{1}$, tiempo en el que el sistema salta a un nuevo estado $x_{1}\neq x_{0}$. Puede ocurrir que el sistema permanezca en $x_{0}$ de manera indefinida, en este caso hacemos $\tau_{1}=\infty$. Si $\tau_{1}$ es finito, el sistema permanecer\'a en $x_{1}$ hasta $\tau_{2}$, y as\'i sucesivamente.
Sea
\begin{equation}
X\left(t\right)=\left\{\begin{array}{cc}
x_{0} & 0\leq t<\tau_{1}\\
x_{1} & \tau_{1}\leq t<\tau_{2}\\
x_{2} & \tau_{2}\leq t<\tau_{3}\\
\vdots &\\
\end{array}\right.
\end{equation}

A este proceso  se le llama {\em proceso de salto}. \index{Proceso de Salto}Si
\begin{equation}
lim_{n\rightarrow\infty}\tau_{n}=\left\{\begin{array}{cc}
<\infty & X_{t}\textrm{ explota,}\\
=\infty & X_{t}\textrm{ no explota.}\\
\end{array}\right.
\end{equation}

Un proceso puro de saltos es un proceso de saltos que satisface la propiedad de Markov.

\begin{Prop}
Un proceso de saltos es Markoviano si y s\'olo si todos los estados no absorbentes $x$ son tales que
\begin{eqnarray*}
P_{x}\left(\tau_{1}>t+s|\tau_{1}>s\right)=P_{x}\left(\tau_{1}>t\right),
\end{eqnarray*}
para $s,t\geq0$, equivalentemente,

\begin{equation}\label{Eq.5}
\frac{1-F_{x}\left(t+s\right)}{1-F_{x}\left(s\right)}=1-F_{x}\left(t\right).
\end{equation}
\end{Prop}

\begin{Note}
Una distribuci\'on $F_{x}$ satisface la ecuaci\'on (\ref{Eq.5}) si y s\'olo si es una funci\'on de distribuci\'on exponencial para todos los estados no absorbentes $x$.
\end{Note}

Por un proceso de nacimiento y muerte \index{Proceso de Nacimiento y Muerte} se entiende un proceso de Markov de Saltos, $\left\{X_{t}\right\}_{t\geq0}$ en $E=\nat$, tal que del estado $n$ s\'olo se puede mover a $n-1$ o $n+1$, es decir, la matriz intensidad \index{Matriz Intensidad}es de la forma:

\begin{equation}
\Lambda=\left(\begin{array}{ccccc}
-\beta_{0}&\beta_{0} & 0 & 0 & \ldots\\
\delta_{1}&-\beta_{1}-\delta_{1} & \beta_{1}&0&\ldots\\
0&\delta_{2}&-\beta_{2}-\delta_{2} & \beta_{2}&\ldots\\
\vdots & & & \ddots &
\end{array}\right)
\end{equation}

donde $\beta_{n}$ son las probabilidades de nacimiento y $\delta_{n}$ las probabilidades de muerte.

La matriz de transici\'on es
\begin{equation}
Q=\left(\begin{array}{ccccc}
0& 1 & 0 & 0 & \ldots\\
q_{1}&0 & p_{1}&0&\ldots\\
0&q_{2}&0& p_{2}&\ldots\\
\vdots & & & \ddots &
\end{array}\right)
\end{equation}
con 
\begin{eqnarray}
\begin{array}{ll}
p_{n}=\frac{\beta_{n}}{\beta_{n}+\delta_{n}}\textrm{  y}& q_{n}=\frac{\delta_{n}}{\beta_{n}+\delta_{n}}.
\end{array}
\end{eqnarray}

\begin{Prop}\label{Prop.2.1}
La recurrencia de un Proceso Markoviano de Saltos $\left\{X_{t}\right\}_{t\geq0}$ con espacio de estados numerable, o equivalentemente de la cadena encajada $\left\{Y_{n}\right\}$ es equivalente a
\begin{equation}\label{Eq.2.1}
\sum_{n=1}^{\infty}\frac{\delta_{1}\cdots\delta_{n}}{\beta_{1}\cdots\beta_{n}}=\sum_{n=1}^{\infty}\frac{q_{1}\cdots q_{n}}{p_{1}\cdots p_{n}}=\infty.
\end{equation}
\end{Prop}

\begin{Lema}\label{Lema.2.2}
Independientemente de la recurrencia o transitoriedad de la cadena, hay una y s\'olo una, salvo m\'ultiplos, soluci\'on $\nu$, a $\nu\Lambda=0$, dada por
\begin{equation}\label{Eq.2.2}
\nu_{n}=\frac{\beta_{0}\cdots\beta_{n-1}}{\delta_{1}\cdots\delta_{n}}\nu_{0}.
\end{equation}
\end{Lema}

\begin{Cor}\label{Corolario2.3}
En el caso recurrente, la medida estacionaria $\mu$ para $\left\{Y_{n}\right\}$, est\'a dada por
\begin{equation}\label{Eq.2.3}
\mu_{n}=\frac{p_{1}\cdots p_{n-1}}{q_{1}\cdots q_{n}}\mu_{0}\textrm{, para }n=1,2,\ldots.
\end{equation}
\end{Cor}

Se define a $S=1+\sum_{n=1}^{\infty}\frac{\beta_{0}\beta_{1}\cdots\beta_{n-1}}{\delta_{1}\delta_{2}\cdots\delta_{n}}$

\begin{Cor}\label{Cor.2.4}
$\left\{X_{t}\right\}$ es erg\'odica si y s\'olo si la ecuaci\'on (\ref{Eq.2.1}) se cumple y adem\'as $S<\infty$, en cuyo caso la distribuci\'on erg\'odica, $\pi$, est\'a dada por
\begin{equation}\label{Eq.2.4}
\pi_{0}=\frac{1}{S}\textrm{, }\pi_{n}=\frac{1}{S}\frac{\beta_{0}\cdots\beta_{n-1}}{\delta_{1}\cdots\delta_{n}}
\end{equation}
para $n=1,2,\ldots$.
\end{Cor}

\begin{Def}
Un proceso irreducible recurrente con medida estacionaria con masa finita es llamado erg\'odico.
\end{Def}

\begin{Teo}\label{Teo4.3}
Un Proceso de Saltos de Markov irreducible no explosivo es erg\'odico si y s\'olo si uno puede encontrar una soluci\'on $\pi$ de probabilidad, $|\pi|=1$, $0\leq\pi_{j}\leq1$, para $\nu\Lambda=0$. En este caso $\pi$ es la distribuci\'on estacionaria.\index{Distribuci\'on Estacionaria}
\end{Teo}
\begin{Cor}\label{Corolario2.4}\index{Cadena Erg\'odica}
$\left\{X_{t}\right\}_{t\geq0}$ es erg\'odica si y s\'olo si (\ref{Eq.2.1}) se cumple y $S<\infty$, en cuyo caso la distribuci\'on estacionaria $\pi$ est\'a dada por

\begin{equation}\label{Eq.2.4}
\pi_{0}=\frac{1}{S}\textrm{,
}\pi_{n}=\frac{1}{S}\frac{\beta_{0}\cdots\beta_{n-1}}{\delta_{1}\cdots\delta_{n}}\textrm{,
}n=1,2,\ldots
\end{equation}
\end{Cor}


Sea $E$ espacio discreto de estados, finito o numerable, y sea $\left\{X_{t}\right\}$ un proceso de Markov con espacio de estados $E$y sea una medida $\mu$ en $E$ definida por sus probabilidades puntuales $\mu_{i}$, escribimos $p_{ij}^{t}=P^{t}\left(i,\left\{j\right\}\right)=P_{i}\left(X_{t}=j\right)$.\\

El monto del tiempo gastado en cada estado es positivo, de modo tal que las trayectorias muestrales son constantes por partes. Para un proceso de saltos denotamos por los tiempos de saltos a $S_{0}=0<S_{1}<S_{2}\cdots$, los tiempos entre saltos consecutivos $T_{n}=S_{n+1}-S_{n}$ y la secuencia de estados visitados por $Y_{0},Y_{1},\ldots$, as\'i las trayectorias muestrales son constantes entre $S_{n}$ consecutivos, continua por la derecha, es decir, $X_{S_{n}}=Y_{n}$.  La descripci\'on de un modelo pr\'actico est\'a dado usualmente en t\'erminos de las intensidades $\lambda\left(i\right)$ y las probabilidades de salto $q_{ij}$ m\'as que en t\'erminos de la matriz de transici\'on $P^{t}$. Sup\'ongase de ahora en adelante que $q_{ii}=0$ cuando $\lambda\left(i\right)>0$

\begin{Teo}
Cualquier Proceso de Markov de Saltos satisface la Propiedad Fuerte de Markov
\end{Teo}

\begin{Teo}\label{Teo.4.2}
Supongamos que $\left\{X_{t}\right\}$ es irreducible recurrente en $E$. Entonces existe una y s\'olo una, salvo m\'ultiplos, medida estacionaria $v$. Esta $v$ tiene la propiedad de que $0<v_{j}<\infty$ para todo $j$ y puede encontrarse en cualquiera de las siguientes formas

\begin{itemize}
\item[i)] Para alg\'un estado $i$, fijo pero arbitrario, $v_{j}$ es el tiempo esperado utilizado en $j$ entre dos llegadas consecutivas al estado $i$;
\begin{equation}\label{Eq.4.2}
v_{j}=\esp_{i}\int_{0}^{w\left(i\right)}\indora\left(X_{t}=j\right)dt
\end{equation}
con $w\left(i\right)=\inf\left\{t>0:X_{t}=i,X_{t^{-}}=\lim_{s\uparrow t}X_{s}\neq i\right\}$. 
\item[ii)]
$v_{j}=\frac{\mu_{j}}{\lambda\left(j\right)}$, donde $\mu$ es estacionaria para $\left\{Y_{n}\right\}$. \item[iii)] como
soluci\'on de $v\Lambda=0$.
\end{itemize}
\end{Teo}

\begin{Def}\index{Proceso Erg\'odico}
Un proceso irreducible recurrente con medida estacionaria de masa finita es llamado erg\'odico.
\end{Def}

\begin{Teo}\label{Teo.4.3}
Un proceso de Markov de saltos irreducible no explosivo es erg\'odico si y s\'olo si se puede encontrar una soluci\'on, de probabilidad, $\pi$, con $|\pi|=1$ y $0\leq\pi_{j}\leq1$, a $\pi\Lambda=0$. En este caso $\pi$ es la distribuci\'on estacionaria.
\end{Teo}

\begin{Cor}\label{Cor.4.4}
Una condici\'on suficiente para la ergodicidad de un proceso irreducible es la existencia de una probabilidad $\pi$ que resuelva el sistema $\pi\Lambda=0$ y que adem\'as tenga la propiedad de que $\sum\pi_{j}\lambda\left(j\right)<\infty$.
\end{Cor}

\begin{Def}\index{Matriz Intensidad}
La matriz intensidad $\Lambda=\left(\lambda\left(i,j\right)\right)_{i,j\in E}$ del proceso de saltos $\left\{X_{t}\right\}_{t\geq0}$ est\'a dada por
\begin{eqnarray}
\begin{array}{l}
\lambda\left(i,j\right)=\lambda\left(i\right)q_{i,j}\textrm{,    }j\neq i\\
\lambda\left(i,i\right)=-\lambda\left(i\right)
\end{array}
\end{eqnarray}
\end{Def}

\begin{Prop}\label{Prop.3.1}
Una matriz $E\times E$, $\Lambda$ es la matriz de intensidad de un proceso markoviano de saltos $\left\{X_{t}\right\}_{t\geq0}$ si y s\'olo si
\begin{eqnarray}
\lambda\left(i,i\right)\leq0\textrm{, }\lambda\left(i,j\right)\textrm{,   }i\neq j\textrm{,  }\sum_{j\in E}\lambda\left(i,j\right)=0.
\end{eqnarray}
\end{Prop}

Para el caso particular de la Cola $M/M/1$, la matr\'iz de intensidad est\'a dada por
\begin{eqnarray*}
\Lambda=\left[\begin{array}{cccccc}
-\beta & \beta & 0 &0 &0& \cdots\\
\delta & -\beta-\delta & \beta & 0 & 0 &\cdots\\
0 & \delta & -\beta-\delta & \beta & 0 &\cdots\\
\vdots & & & & & \ddots\\
\end{array}\right]
\end{eqnarray*}

\begin{Prop}
Si el proceso es erg\'odico, entonces existe una versi\'on estrictamente estacionaria $\left\{X_{t}\right\}_{-\infty<t<\infty}$con doble tiempo infinito.
\end{Prop}

\begin{Teo}
Si $\left\{X_{t}\right\}$ es erg\'odico y $\pi$ es la distribuci\'on estacionaria, entonces para todo $i,j$, $p_{ij}^{t}\rightarrow\pi_{j}$ cuando $t\rightarrow\infty$.
\end{Teo}

\begin{Cor}
Si $\left\{X_{t}\right\}$ es irreducible recurente pero no erg\'odica, es decir $|v|=\infty$, entonces $p_{ij}^{t}\rightarrow0$ para todo $i,j\in E$.
\end{Cor}

\begin{Cor}
Para cualquier proceso Markoviano de Saltos minimal, irreducible o no, los l\'imites $li_{t\rightarrow\infty}p_{ij}^{t}$ existe.
\end{Cor}

\section{Notaci\'on Kendall-Lee}
Dado un sistema de espera (colas) a partir de este momento se har\'an las siguientes consideraciones:
\begin{itemize}
\item[a) ]Si $t_{n}$ es el tiempo aleatorio en el que llega al sistema el $n$-\'esimo cliente, para $n=1,2,\ldots$, $t_{0}=0$ y $t_{0}<t_{1}<\cdots$ se definen los tiempos entre arribos $\tau_{n}=t_{n}-t_{n-1}$ para $n=1,2,\ldots$, variables aleatorias independientes e id\'enticamente distribuidas.

\item[b) ]Los tiempos entre arribos tienen un valor medio $0<E\left(\tau\right)=\frac{1}{\beta}<\infty$, es decir, $\beta$ se puede ver como la tasa o intensidad promedio de arribos al sistema por unidad de tiempo.
\item[c) ]  Adem\'as se supondr\'a que los servidores son identicos y si $s$ denota la variable aleatoria que describe el tiempo de servicio, entonces $E\left(s\right)=\frac{1}{\delta}$, $\delta$ es la tasa promedio de servicio por servidor.
\end{itemize}

La notaci\'on de Kendall-Lee es una forma abreviada de describir un sistema de espera con las siguientes componentes:
\begin{itemize}
\item[a)] {\em\bf Fuente}: Poblaci\'on de clientes potenciales del sistema, esta puede ser finita o infinita. 
\item[b)] {\em\bf Proceso de Arribos}: Proceso determinado por la funci\'on de distribuci\'on $A\left(t\right)=P\left\{\tau\leq t\right\}$ de los tiempos entre arribos.
\end{itemize}

Adem\'as tenemos las siguientes igualdades
\begin{equation}\label{Eq.0.1}
N\left(t\right)=N_{q}\left(t\right)+N_{s}\left(s\right)
\end{equation}
donde
\begin{itemize}
\item[a) ] $N\left(t\right)$ es el n\'umero de clientes en el sistema al tiempo $t$. 
\item[b) ] $N_{q}\left(t\right)$ es el n\'umero de cliente en la cola al tiempo $t$.
\item[c) ] $N_{s}\left(t\right)$ es el n\'umero de clientes recibiendo servicio en el tiempo $t$.
\end{itemize}

Bajo la hip\'otesis de estacionareidad, es decir, las caracter\'isticas de funcionamiento del sistema se han estabilizado en valores independientes del tiempo, entonces
\begin{equation}
N=N_{q}+N_{s}.
\end{equation}

Los valores medios de las cantidades anteriores se escriben como 
\begin{eqnarray}
\begin{array}{ccc}
L=E\left(N\right), &L_{q}=E\left(N_{q}\right)\textrm{ y }&L_{s}=E\left(N_{s}\right),
\end{array}
\end{eqnarray}
entonces de la ecuaci\'on \ref{Eq.0.1} se obtiene

\begin{equation}
L=L_{q}+L_{s}
\end{equation}
Si $q$ es el tiempo que pasa un cliente en la cola antes de recibir servicio, y W es el tiempo total que un cliente pasa en el sistema, entonces \begin{eqnarray*}w=q+s\end{eqnarray*} por lo tanto 
\begin{equation}W=W_{q}+W_{s},\end{equation} donde \begin{eqnarray*}W=E\left(w\right), W_{q}=E\left(q\right)\textrm{ y }W_{s}=E\left(s\right)=\frac{1}{\delta}\end{eqnarray*}.

La intensidad de tr\'afico se define como
\begin{equation}
\rho=\frac{E\left(s\right)}{E\left(\tau\right)}=\frac{\beta}{\delta}.
\end{equation}

La utilizaci\'on por servidor es
\begin{eqnarray}
u=\frac{\rho}{c}=\frac{\beta}{c\delta}\textrm{, donde }c\textrm{ es el n\'umero de servidores.}
\end{eqnarray}

La siguiente notaci\'on es una forma abreviada de describir un sistema de espera con componentes dados:

\begin{equation}\label{Notacion.K.L.}
A/S/c/K/F/d
\end{equation}

Cada una de las letras describe:

\begin{itemize}
\item $A$ es la distribuci\'on de los tiempos entre arribos.
\item $S$ es la distribuci\'on del tiempo de servicio.
\item $c$ es el n\'umero de servidores.
\item $K$ es la capacidad del sistema.
\item $F$ es el n\'umero de individuos en la fuente.
\item $d$ es la disciplina del servicio
\end{itemize}

Usualmente se acostumbra suponer que la capacidad del sistema es infinita, $K=\infty$, la cantidad de personasen espera es infinita $F=\infty$, y la pol\'itica de servicio es de tipo $d=FIFO$, es decir, \textit{First In First Out}. Las distribuciones usuales para $A$ y $B$ son:

\begin{itemize}
\item $GI$ para la distribuci\'on general de los tiempos entre arribos.
\item $G$ distribuci\'on general del tiempo de servicio.
\item $M$ Distribuci\'on exponencial para $A$ o $S$.
\item $E_{K}$ Distribuci\'on Erlang-$K$, para $A$ o $S$.
\item $D$ tiempos entre arribos o de servicio constantes, es decir, deterministicos.
\end{itemize}

\subsection*{Cola $M/M/1$}\index{Cola $M/M/1$}
Este modelo corresponde a un proceso de nacimiento y muerte con $\beta_{n}=\beta$ y $\delta_{n}=\delta$ independiente del valor de $n$. La intensidad de tr\'afico $\rho=\frac{\beta}{\delta}$, implica que el criterio de recurrencia (ecuaci\'on \ref{Eq.2.1}) quede de la forma:
\begin{eqnarray}
1+\sum_{n=1}^{\infty}\rho^{-n}=\infty.
\end{eqnarray}
Equivalentemente el proceso es recurrente si y s\'olo si
\begin{eqnarray}
\sum_{n\geq1}\left(\frac{\beta}{\delta}\right)^{n}<\infty\Leftrightarrow \frac{\beta}{\delta}<1.
\end{eqnarray}
Entonces
$S=\frac{\delta}{\delta-\beta}$, luego por la ecuaci\'on \ref{Eq.2.4} se tiene que
\begin{eqnarray}
\begin{array}{ll}
\pi_{0}=\frac{\delta-\beta}{\delta}=1-\frac{\beta}{\delta},\\
\pi_{n}=\pi_{0}\left(\frac{\beta}{\delta}\right)^{n}=\left(1-\frac{\beta}{\delta}\right)\left(\frac{\beta}{\delta}\right)^{n}=\left(1-\rho\right)\rho^{n}.
\end{array}
\end{eqnarray}

Lo cual nos lleva a la siguiente proposici\'on:

\begin{Prop}
La cola $M/M/1$ con intensidad de tr\'afico $\rho$, es recurrente si y s\'olo si $\rho\leq1$.
\end{Prop}

Adem\'as se tiene

\begin{Prop}
La cola $M/M/1$ con intensidad de tr\'afico $\rho$ es erg\'odica si y s\'olo si $\rho<1$. En cuyo caso, la distribuci\'on de equilibrio $\pi$ de la longitud de la cola es geom\'etrica, 
\begin{eqnarray}
pi_{n}=\left(1-\rho\right)\rho^{n}\textrm{, para }n=1,2,\ldots.
\end{eqnarray}
\end{Prop}

De la proposici\'on anterior se desprenden varios hechos importantes.
\begin{itemize}
\item[a) ] $\prob\left[X_{t}=0\right]=\pi_{0}=1-\rho$, es decir, la probabilidad de que el sistema se encuentre ocupado.
\item[b) ] De las propiedades de la distribuci\'on Geom\'etrica se desprende que
\begin{eqnarray}
\begin{array}{l}
\esp\left[X_{t}\right]=\frac{\rho}{1-\rho},\\
Var\left[X_{t}\right]=\frac{\rho}{\left(1-\rho\right)^{2}}.
\end{array}
\end{eqnarray}
\end{itemize}

Si $L$ es el n\'umero esperado de clientes en el sistema, incluyendo los que est\'an siendo atendidos, entonces
\begin{eqnarray}
L=\frac{\rho}{1-\rho}.
\end{eqnarray}
Si adem\'as $W$ es el tiempo total del cliente en la cola:
\begin{eqnarray}
W=W_{q}+W_{s},\\
\rho=\frac{\esp\left[s\right]}{\esp\left[\tau\right]}=\beta W_{s},
\end{eqnarray}
puesto que $W_{s}=\esp\left[s\right]$ y $\esp\left[\tau\right]=\frac{1}{\delta}$. Por la f\'ormula de Little
\begin{equation}
L=\lambda W,
\end{equation}
\begin{eqnarray}
W=\frac{L}{\beta}=\frac{\frac{\rho}{1-\rho}}{\beta}=\frac{\rho}{\delta}\frac{1}{1-\rho}=\frac{W_{s}}{1-\rho}=\frac{1}{\delta\left(1-\rho\right)}=\frac{1}{\delta-\beta},
\end{eqnarray}

luego entonces

\begin{eqnarray}
W_{q}=W-W_{s}=\frac{1}{\delta-\beta}-\frac{1}{\delta}=\frac{\beta}{\delta(\delta-\beta)}=\frac{\rho}{1-\rho}\frac{1}{\delta}=\esp\left[s\right]\frac{\rho}{1-\rho}.
\end{eqnarray}

Entonces

\begin{eqnarray}
L_{q}=\beta W_{q}=\frac{\rho^{2}}{1-\rho}.
\end{eqnarray}

Finalmente, tenemos las siguientes proposiciones:

\begin{Prop}
\begin{eqnarray}
\begin{array}{l}
W\left(t\right)=1-e^{-\frac{t}{W}}.\\
W_{q}\left(t\right)=1-\rho\exp^{-\frac{t}{W}}\textrm{, donde }W=\esp(w).
\end{array}
\end{eqnarray}

\end{Prop}

\begin{Prop}
La cola M/M/1 con intensidad de tr\'afico $\rho$ es recurrente si y s\'olo si $\rho\leq1$
\end{Prop}

\begin{Prop}
La cola M/M/1 con intensidad de tr\'afica $\rho$ es ergodica si y s\'olo si $\rho<1$. En este caso, la distribuci\'on de equilibrio $\pi$ de la longitud de la cola es geom\'etrica, 
\begin{eqnarray}
\pi_{n}=\left(1-\rho\right)\rho^{n}\textrm{, para }n=0,1,2,\ldots.
\end{eqnarray}
\end{Prop}
\subsection*{Cola $M/M/\infty$}\index{Cola $M/M/\infty$}

Este tipo de modelos se utilizan para estimar el n\'umero de l\'ineas en uso en una gran red comunicaci\'on o para estimar valores en los sistemas $M/M/c$ o $M/M/c/c$, se puede pensar que siempre hay un servidor disponible para cada cliente que llega.

Se puede considerar como un proceso de nacimiento y muerte con par\'ametros $\beta_{n}=\beta$ y $\mu_{n}=n\mu$ para $n=0,1,2,\ldots$. Este modelo corresponde a $\beta_{n}=\beta$ y $\delta_{n}=n\delta$, en este caso, el par\'ametro de inter\'es $\eta=\frac{\beta}{\delta}$, luego, la ecuaci\'on (\ref{Eq.2.1}) queda de la forma:

\begin{eqnarray}
\sum_{n=1}^{\infty}\frac{\delta_{1}\cdots\delta_{n}}{\beta_{1}\cdots\beta_{n}}=\sum_{n=1}^{\infty}n!\eta^{-n}=\infty\textrm{con }S=1+\sum_{n=1}^{\infty}\frac{\eta^{n}}{n!}=e,
\end{eqnarray}
 entonces por la ecuaci\'on (\ref{Eq.2.4}) se tiene que

\begin{eqnarray}\label{MMinf.pi}
\begin{array}{l}
\pi_{0}=e^{\rho},\\
\pi_{n}=e^{-\rho}\frac{\rho^{n}}{n!}.
\end{array}
\end{eqnarray}
Entonces, el n\'umero promedio de servidores ocupados es equivalente a considerar el n\'umero de clientes en el  sistema, es decir,
\begin{eqnarray}
\begin{array}{l}
L=\esp\left[N\right]=\rho,\\
Var\left[N\right]=\rho.
\end{array}
\end{eqnarray}
Adem\'as se tiene que $W_{q}=0$ y $L_{q}=0$. El tiempo promedio en el sistema es el tiempo promedio de servicio, es decir, 
\begin{equation}
W=\esp\left[s\right]=\frac{1}{\delta}.
\end{equation}
Resumiendo, tenemos la sisuguiente proposici\'on:

\begin{Prop}
La cola $M/M/\infty$ es erg\'odica para todos los valores de $\eta$. La distribuci\'on de equilibrio $\pi$ es Poisson con media $\eta$,
\begin{eqnarray}
\pi_{n}=\frac{e^{-n}\eta^{n}}{n!}.
\end{eqnarray}
\end{Prop}
\subsection*{Cola $M/M/m$}\index{Cola $M/M/m$}

Este sistema considera $m$ servidores id\'enticos, con tiempos entre arribos y de servicio exponenciales con medias $\esp\left[\tau\right]=\frac{1}{\beta}$ y $\esp\left[s\right]=\frac{1}{\delta}$. definimos ahora la utilizaci\'on por servidor como $u=\frac{\rho}{m}$ que tambi\'en se puede interpretar como la fracci\'on de tiempo promedio que cada servidor est\'a ocupado. La cola $M/M/m$ se puede considerar como un proceso de nacimiento y muerte con par\'ametros: $\beta_{n}=\beta$ para $n=0,1,2,\ldots$ y
\begin{eqnarray}
\delta_{n}=\left\{\begin{array}{cc}
n\delta & n=0,1,\ldots,m-1\\
c\delta & n=m,\ldots\\
\end{array}\right.
\end{eqnarray}

entonces  la condici\'on de recurrencia se va a cumplir s\'i y s\'olo si 
\begin{eqnarray}
\sum_{n\geq1}\frac{\beta_{0}\cdots\beta_{n-1}}{\delta_{1}\cdots\delta_{n}}<\infty,
\end{eqnarray}
equivalentemente se debe de cumplir que
\begin{eqnarray}
S=1+\sum_{n\geq1}\frac{\beta_{0}\cdots\beta_{n-1}}{\delta_{1}\cdots\delta_{n}}=\sum_{n=0}^{m-1}\frac{\beta_{0}\cdots\beta_{n-1}}{\delta_{1}\cdots\delta_{n}}+\sum_{n=0}^{\infty}\frac{\beta_{0}\cdots\beta_{n-1}}{\delta_{1}\cdots\delta_{n}}=\sum_{n=0}^{m-1}\frac{\beta^{n}}{n!\delta^{n}}+\sum_{n=0}^{\infty}\frac{\rho^{m}}{m!}u^{n}
\end{eqnarray}
converja, lo cual ocurre si $u<1$, en este caso

\begin{eqnarray}
S=\sum_{n=0}^{m-1}\frac{\rho^{n}}{n!}+\frac{\rho^{m}}{m!}\left(1-u\right)
\end{eqnarray}
luego, para este caso se tiene que

\begin{eqnarray}
\begin{array}{l}
\pi_{0}=\frac{1}{S}\\
\pi_{n}=\left\{\begin{array}{cc}
\pi_{0}\frac{\rho^{n}}{n!} & n=0,1,\ldots,m-1\\
\pi_{0}\frac{\rho^{n}}{m!m^{n-m}}& n=m,\ldots\\
\end{array}\right.

\end{array}
\end{eqnarray}
Al igual que se hizo antes, determinaremos los valores de
$L_{q},W_{q},W$ y $L$:

\begin{eqnarray}
\begin{array}{l}
L_{q}=\esp\left[N_{q}\right]=\sum_{n=0}^{\infty}\left(n-m\right)\pi_{n}=\sum_{n=0}^{\infty}n\pi_{n+m}=\sum_{n=0}^{\infty}n\pi_{0}\frac{\rho^{n+m}}{m!m^{n+m}}=\pi_{0}\frac{\rho^{m}}{m!}\sum_{n=0}^{\infty}nu^{n}\\
=\pi_{0}\frac{u\rho^{m}}{m!}\sum_{n=0}^{\infty}\frac{d}{du}u^{n}=\pi_{0}\frac{u\rho^{m}}{m!}\frac{d}{du}\sum_{n=0}^{\infty}u^{n}=\pi_{0}\frac{u\rho^{m}}{m!}\frac{d}{du}\left(\frac{1}{1-u}\right)=\pi_{0}\frac{u\rho^{m}}{m!}\frac{1}{\left(1-u\right)^{2}},
\end{array}
\end{eqnarray}

es decir
\begin{equation}
L_{q}=\frac{u\pi_{0}\rho^{m}}{m!\left(1-u\right)^{2}},
\end{equation}
luego
\begin{equation}
W_{q}=\frac{L_{q}}{\beta}.
\end{equation}
Adem\'as
\begin{equation}
W=W_{q}+\frac{1}{\delta}
\end{equation}

Si definimos
\begin{eqnarray}
C\left(m,\rho\right)=\frac{\pi_{0}\rho^{m}}{m!\left(1-u\right)}=\frac{\pi_{m}}{1-u},
\end{eqnarray}
que es la probabilidad de que un cliente que llegue al sistema
tenga que esperar en la cola. Entonces podemos reescribir las
ecuaciones reci\'en enunciadas:

\begin{eqnarray}
\begin{array}{ll}
L_{q}=\frac{C\left(m,\rho\right)u}{1-u},&\textrm{y }
W_{q}=\frac{C\left(m,\rho\right)\esp\left[s\right]}{m\left(1-u\right)}
\end{array}
\end{eqnarray}
Por tanto tenemos las siguientes proposiciones:

\begin{Prop}
La cola $M/M/m$ con intensidad de tr\'afico $\rho$ es erg\'odica si y s\'olo si $\rho<1$. En este caso la distribuci\'on erg\'odica $\pi$ est\'a dada por
\begin{eqnarray}
\pi_{n}=\left\{\begin{array}{cc}
\frac{1}{S}\frac{\eta^{n}}{n!} & 0\leq n\leq m,\\
\frac{1}{S}\frac{\eta^{m}}{m!}\rho^{n-m} & m\leq n<\infty.
\end{array}\right.
\end{eqnarray}
\end{Prop}

\begin{Prop}
Para $t\geq0$
\begin{itemize}
\item[a)]
\begin{eqnarray}
W_{q}\left(t\right)=1-C\left(m,\rho\right)e^{-c\delta
t\left(1-u\right)}.
\end{eqnarray} 
\item[b)]\begin{eqnarray}
W\left(t\right)=\left\{\begin{array}{cc}
1+e^{-\delta t}\frac{\rho-m+W_{q}\left(0\right)}{m-1-\rho}+e^{-m\delta t\left(1-u\right)}\frac{C\left(m,\rho\right)}{m-1-\rho} & \rho\neq m-1,\\
1-\left(1+C\left(m,\rho\right)\delta t\right)e^{-\delta t} & \rho=m-1.
\end{array}\right.
\end{eqnarray}
\end{itemize}
\end{Prop}

Resumiendo, para este caso $\beta_{n}=\beta$ y $\delta_{n}=m\left(n\right)\delta$, donde $m\left(n\right)$ es el n\'umero de servidores ocupados en el estado $n$, es decir $m\left(n\right)=m$, para $n\geq m$ y $m\left(n\right)=m$ para $1\leq n\leq m$. La intensidad de tr\'afico es 
\begin{eqnarray}
\rho=\frac{\beta}{m\delta}\textrm{ y }\frac{\beta_{n}}{\delta_{n}}=\rho,\textrm{ para }n\geq m.
\end{eqnarray}
As\'i, al igual que en el caso $m=1$, la ecuaci\'on (\ref{Eq.2.1}) y la recurrencia se cumplen si y s\'olo si 
\begin{eqnarray}
\sum_{n=1}^{\infty}\rho^{-n}=\infty\textrm{, es decir, cuando }\rho\leq1.
\end{eqnarray}
 
\subsection*{Cola $M/M/m/m$}\index{Cola $M/M/m/m$}

Consideremos un sistema con $m$ servidores id\'enticos, pero ahora cada uno es de capacidad finita $m$. Si todos los servidores se encuentran ocupados, el siguiente usuario en llegar se pierde pues no se le deja esperar a que reciba servicio. Este tipo de sistemas pueden verse como un proceso de nacimiento y muerte con
\begin{eqnarray}
\begin{array}{ll}
\beta_{n}=\left\{\begin{array}{cc}
\beta & n=0,1,2,\ldots,m-1,\\
0 & n\geq m.\\
\end{array}
\right.,&
\delta_{n}=\left\{\begin{array}{cc}
n\delta & n=0,1,2,\ldots,m-1,\\
0 & n\geq m.\\
\end{array}
\right.
\end{array}
\end{eqnarray}

El proceso tiene espacio de estados finitos, $S=\left\{0,1,\ldots,m\right\}$, entonces de las ecuaciones que determinan la distribuci\'on estacionaria se tiene que
\begin{eqnarray}\label{Eq.13.1}
\begin{array}{ll}
\pi_{n}=\left\{\begin{array}{cc}
\pi_{0}\frac{\rho^{n}}{n!} & n=0,1,2,\ldots,m\\
0 & n\geq m\\
\end{array}
\right.,&
\textrm{y adem\'as, }
\pi_{0}=\left(\sum_{n=0}^{m}\frac{\rho^{n}}{n!}\right)^{-1}.
\end{array}
\end{eqnarray}
A la ecuaci\'on (\ref{Eq.13.1}) se le llama {\em distribuci\'on truncada}. Si definimos 
\begin{eqnarray}
\pi_{m}=B\left(m,\rho\right)=\pi_{0}\frac{\rho^{m}}{m!},
\end{eqnarray}
$\pi_{m}$ representa la probabilidad de que todos los servidores se encuentren ocupados, y tambi\'en se le conoce como {\em f\'ormula de p\'erdida de Erlang}. Necesariamente en este caso el tiempo de espera en la cola $W_{q}$ y el n\'umero promedio de clientes en la cola $L_{q}$, deben de ser cero puesto que no se permite esperar para recibir servicio, m\'as a\'un, el tiempo de espera en el sistema y el tiempo de serivcio tienen la misma distribuci\'on, es decir,
\begin{eqnarray}
W\left(t\right)=\prob\left\{w\leq t\right\}=1-e^{-\mu t},
\end{eqnarray}
en particular
\begin{eqnarray}
W=\esp\left[w\right]=\esp\left[s\right]=\frac{1}{\delta}.
\end{eqnarray}
Por otra parte, el n\'umero esperado de clientes en el sistema es
\begin{eqnarray}
L=\esp\left[N\right]=\sum_{n=0}^{m}n\pi_{n}=\pi_{0}\rho\sum_{n=0}^{m}\frac{\rho^{n-1}}{\left(n-1\right)!}=\pi_{0}\rho\sum_{n=0}^{m-1}\frac{\rho^{n}}{n!},
\end{eqnarray}
entonces, se tiene que
\begin{equation}
L=\rho\left(1-B\left(m,\rho\right)\right)=\esp\left[s\right]\left(1-B\left(m,\rho\right)\right).
\end{equation}
Adem\'as
\begin{equation}
\delta_{q}=\delta\left(1-B\left(m,\rho\right)\right),
\end{equation}
que representa la tasa promedio efectiva de arribos al sistema.
\subsection*{Cola $M/G/1$}\index{Cola $M/G/1$}

Consideremos un sistema de espera con un servidor, en el que los tiempos entre arribos son exponenciales, y los tiempos de servicio tienen una distribuci\'on general $G$. Sea $N\left(t\right)$, con $t\geq0$, el n\'umero de clientes en el sistema al tiempo $t$, y sean $t_{1}<t_{2}<\dots$ los tiempos sucesivos en los que los clientes completan su servicio y salen del sistema.\\

La sucesi\'on $\left\{X_{n}\right\}$ definida por
$X_{n}=N\left(t_{n}\right)$ es una cadena de Markov, en espec\'ifico es la cadena encajada\index{Cadena Encajada} del proceso de llegadas de usuarios. Sea $U_{n}$ el n\'umero de clientes que llegan al sistema durante el tiempo de servicio del $n$-\'esimo cliente, entonces se tiene que

\begin{eqnarray}
X_{n+1}=\left\{\begin{array}{cc}
X_{n}-1+U_{n+1} & \textrm{si }X_{n}\geq1,\\
U_{n+1} & \textrm{si }X_{n}=0.
\end{array}\right.
\end{eqnarray}

Dado que los procesos de arribos de los usuarios es Poisson con par\'ametro $\lambda$\index{Proceso Poisson}, la probabilidad condicional de que lleguen $j$ clientes al sistema dado que el tiempo de servicio es $s=t$, resulta:
\begin{eqnarray}
\prob\left\{U=j|s=t\right\}=e^{-\lambda t}\frac{\left(\lambda
t\right)^{j}}{j!}\textrm{,   }j=0,1,\ldots
\end{eqnarray}
por el teorema de la probabilidad total se tiene que
\begin{equation}
a_{j}=\prob\left\{U=j\right\}=\int_{0}^{\infty}\prob\left\{U=j|s=t\right\}dG\left(t\right)=\int_{0}^{\infty}e^{-\lambda
t}\frac{\left(\lambda t\right)^{j}}{j!}dG\left(t\right)
\end{equation}
donde $G$ es la distribuci\'on de los tiempos de servicio. Las probabilidades de transici\'on de la cadena est\'an dadas por
\begin{equation}
p_{0j}=\prob\left\{U_{n+1}=j\right\}=a_{j}\textrm{, para
}j=0,1,\ldots
\end{equation}
y para $i\geq1$
\begin{equation}
p_{ij}=\left\{\begin{array}{cc}
\prob\left\{U_{n+1}=j-i+1\right\}=a_{j-i+1}&\textrm{, para }j\geq i-1\\
0 & j<i-1\\
\end{array}
\right.
\end{equation}
Entonces la matriz de transici\'on es:
\begin{eqnarray*}
P=\left[\begin{array}{ccccc}
a_{0} & a_{1} & a_{2} & a_{3} & \cdots\\
a_{0} & a_{1} & a_{2} & a_{3} & \cdots\\
0 & a_{0} & a_{1} & a_{2} & \cdots\\
0 & 0 & a_{0} & a_{1} & \cdots\\
\vdots & \vdots & \cdots & \ddots &\vdots\\
\end{array}
\right].
\end{eqnarray*}
Sea $\rho=\sum_{n=0}na_{n}$, entonces se tiene el siguiente teorema:
\begin{Teo}
La cadena encajada $\left\{X_{n}\right\}$ es
\begin{itemize}
\item[a)] Recurrente positiva si $\rho<1$,
\item[b)] Transitoria si $\rho>1$, 
\item[c)] Recurrente nula si $\rho=1$.
\end{itemize}
\end{Teo}

Recordemos que si la cadena de Markov $\left\{X_{n}\right\}$ tiene una distribuci\'on estacionaria entonces existe una distribuci\'on de probabilidad $\pi=\left(\pi_{0},\pi_{1},\ldots,\right)$, con $\pi_{i}\geq0$ y $\sum_{i\geq1}\pi_{i}=1$ tal que satisface la
ecuaci\'on $\pi=\pi P$, equivalentemente
\begin{equation}\label{Eq.18.9}
\pi_{j}=\sum_{i=0}^{\infty}\pi_{k}p_{ij},\textrm{ para
}j=0,1,2,\ldots
\end{equation}
que se puede ver como
\begin{equation}\label{Eq.19.6}
\pi_{j}=\pi_{0}a_{j}+\sum_{i=1}^{j+1}\pi_{i}a_{j-i+1}\textrm{,
para }j=0,1,\ldots
\end{equation}
si definimos
\begin{eqnarray}
\pi\left(z\right)=\sum_{j=0}^{\infty}\pi_{j}z^{j}\textrm{, y } 
A\left(z\right)=\sum_{j=0}^{\infty}a_{j}z^{j}\textrm{. con }|z_{j}|\leq1.
\end{eqnarray}
 Si la ecuaci\'on (\ref{Eq.19.6}) la multiplicamos por $z^{j}$ y sumando sobre $j$, se tiene que
\begin{eqnarray*}
\sum_{j=0}^{\infty}\pi_{j}z^{j}&=&\sum_{j=0}^{\infty}\pi_{0}a_{j}z^{j}+\sum_{j=0}^{\infty}\sum_{i=1}^{j+1}\pi_{i}a_{j-i+1}z^{j}=\pi_{0}\sum_{j=0}^{\infty}a_{j}z^{j}+\sum_{j=0}^{\infty}a_{j}z^{j}\sum_{i=1}^{\infty}\pi_{i}a_{i-1}\\
&=&\pi_{0}A\left(z\right)+A\left(z\right)\left(\frac{\pi\left(z\right)-\pi_{0}}{z}\right)
\end{eqnarray*}
es decir,

\begin{equation}
\pi\left(z\right)=\pi_{0}A\left(z\right)+A\left(z\right)\left(\frac{\pi\left(z\right)-\pi_{0}}{z}\right)\Leftrightarrow\pi\left(z\right)=\frac{\pi_{0}A\left(z\right)\left(z-1\right)}{z-A\left(z\right)}
\end{equation}

Si $z\rightarrow1$, entonces $A\left(z\right)\rightarrow A\left(1\right)=1$, y adem\'as $A^{'}\left(z\right)\rightarrow A^{'}\left(1\right)=\rho$. Si aplicamos la Regla de L'Hospital se tiene que
\begin{eqnarray}
\sum_{j=0}^{\infty}\pi_{j}=lim_{z\rightarrow1^{-}}\pi\left(z\right)=\pi_{0}lim_{z\rightarrow1^{-}}\frac{z-1}{z-A\left(z\right)}=\frac{\pi_{0}}{1-\rho}
\end{eqnarray}
Retomando,
\begin{eqnarray*}
a_{j}=\prob\left\{U=j\right\}=\int_{0}^{\infty}e^{-\lambda
t}\frac{\left(\lambda t\right)^{n}}{n!}dG\left(t\right)\textrm{,
para }n=0,1,2,\ldots
\end{eqnarray*}
entonces
\begin{eqnarray}
\begin{array}{l}
\rho=\sum_{n=0}^{\infty}na_{n}=\sum_{n=0}^{\infty}n\int_{0}^{\infty}e^{-\lambda t}\frac{\left(\lambda t\right)^{n}}{n!}dG\left(t\right)=\int_{0}^{\infty}\sum_{n=0}^{\infty}ne^{-\lambda t}\frac{\left(\lambda
t\right)^{n}}{n!}dG\left(t\right)=\int_{0}^{\infty}\lambda tdG\left(t\right)=\lambda\esp\left[s\right].
\end{array}
\end{eqnarray}

Adem\'as, se tiene que $\rho=\beta\esp\left[s\right]=\frac{\beta}{\delta}$, y la distribuci\'on estacionaria est\'a dada por
\begin{eqnarray}
\begin{array}{l}
\pi_{j}=\pi_{0}a_{j}+\sum_{i=1}^{j+1}\pi_{i}a_{j-i+1}\textrm{, para }j=0,1,\ldots\\
\pi_{0}=1-\rho.
\end{array}
\end{eqnarray}
Por otra parte se tiene que\begin{equation}
L=\pi^{'}\left(1\right)=\rho+\frac{A^{''}\left(1\right)}{2\left(1-\rho\right)}
\end{equation}

pero 
\begin{eqnarray*}
A^{''}\left(1\right)&=&\sum_{n=1}n\left(n-1\right)a_{n}= \esp\left[U^{2}\right]-\esp\left[U\right],\\
\esp\left[U\right]&=&\rho,\textrm{ y }\\
\esp\left[U^{2}\right]&=&\lambda^{2}\esp\left[s^{2}\right]+\rho.
\end{eqnarray*}

Por lo tanto 
\begin{equation}
L=\rho+\frac{\beta^{2}\esp\left[s^{2}\right]}{2\left(1-\rho\right)}.
\end{equation}

De las f\'ormulas de Little, se tiene que 
\begin{equation}
W=E\left(w\right)=\frac{L}{\beta},
\end{equation}
tambi\'en el tiempo de espera en la cola
\begin{equation}
W_{q}=\esp\left(q\right)=\esp\left(w\right)-\esp\left(s\right)=\frac{L}{\beta}-\frac{1}{\delta},
\end{equation}
adem\'as el n\'umero promedio de clientes en la cola es
\begin{equation}
L_{q}=\esp\left(N_{q}\right)=\beta W_{q}=L-\rho
\end{equation}

\subsection*{Cola con Infinidad de Servidores}
Este caso corresponde a $\beta_{n}=\beta$ y $\delta_{n}=n\delta$. El par\'ametro de inter\'es es $\eta=\frac{\beta}{\delta}$, de donde se obtiene:
\begin{eqnarray}
\sum_{n\geq0}\frac{\delta_{1}\cdots\delta_{n}}{\beta_{1}\cdots\beta_{n}}&=&\sum_{n=1}^{\infty}n!\eta^{n}=\infty,\\
S&=&1+\sum_{n=1}^{\infty}\frac{\eta^{n}}{n!}=e^{n}.
\end{eqnarray}
Por lo tanto se tiene la siguiente proposici\'on: 
\begin{Prop}
La cola $M/M/\infty$ es ergodica para todos los valores de $\eta$. La distribuci\'on de equilibrio $\pi$ es Poisson con media $\eta$, 
\begin{eqnarray}
pi_{n}=\frac{e^{-n}\eta}{n!}.
\end{eqnarray}
\end{Prop}
\section{Redes de Colas:Sistemas Abiertos}

Consid\'erese un sistema con dos servidores, en los cuales los usuarios llegan de acuerdo a un proceso Poisson con intensidad $\lambda_{1}$ al primer servidor, despu\'es de ser atendido se pasa a la siguiente cola en el segundo servidor. Cada servidor atiende a un usuario a la vez con tiempo exponencial con raz\'on $\mu_{i}$, para $i=1,2$. A este tipo de sistemas se les conoce como sistemas secuenciales.

Def\'inase el par $\left(n,m\right)$ como el n\'umero de usuarios en el servidor 1 y 2 respectivamente. Las ecuaciones de balance son
\begin{eqnarray}\label{Eq.Balance}
\begin{array}{c}
\lambda P_{0,0}=\mu_{2}P_{0,1}\\
\left(\lambda+\mu_{1}\right)P_{n,0}=\mu_{2}P_{n,1}+\lambda P_{n-1,0}\\
\left(\lambda+\mu_{2}\right)P_{0,m}=\mu_{2}P_{0,m+1}+\mu_{1}P_{1,m-1}\\
\left(\lambda+\mu_{1}+\mu_{2}\right)P_{n,m}=\mu_{2}P_{n,m+1}+\mu_{1}P_{n+1,m-1}+\lambda
P_{n-1,m}
\end{array}
\end{eqnarray}

Cada servidor puede ser visto como un modelo de tipo $M/M/1$, de igual manera el proceso de salida de una cola $M/M/1$ con raz\'on $\lambda$, nos permite asumir que el servidor 2 tambi\'en es una cola $M/M/1$. Adem\'as, la probabilidad de que haya $n$ usuarios en el servidor 1 es
\begin{eqnarray}
\begin{array}{l}
P\left\{n\textrm{ en el servidor }1\right\}=\left(\frac{\lambda}{\mu_{1}}\right)^{n}\left(1-\frac{\lambda}{\mu_{1}}\right)=\rho_{1}^{n}\left(1-\rho_{1}\right)\\
P\left\{m\textrm{ en el servidor }2\right\}=\left(\frac{\lambda}{\mu_{2}}\right)^{n}\left(1-\frac{\lambda}{\mu_{2}}\right)=\rho_{2}^{m}\left(1-\rho_{2}\right)
\end{array}
\end{eqnarray}
Si el n\'umero de usuarios en los servidores 1 y 2 son variables aleatorias independientes, se sigue que:
\begin{equation}\label{Eq.8.16}
P_{n,m}=\rho_{1}^{n}\left(1-\rho_{1}\right)\rho_{2}^{m}\left(1-\rho_{2}\right)
\end{equation}
Verifiquemos que $P_{n,m}$ satisface las ecuaciones de balance (\ref{Eq.Balance}) Antes de eso, enunciemos unas igualdades que nos ser\'an de utilidad:

$\mu_{i}\rho_{i}=\lambda$ para $i=1,2.$ Entonces, \[\lambda P_{0,0}=\lambda\left(1-\rho_{1}\right)\left(1-\rho_{2}\right),\] y \[\mu_{2} P_{0,1}=\mu_{2}\left(1-\rho_{1}\right)\rho_{2}\left(1-\rho_{2}\right).\]

Entonces $\lambda P_{0,0}=\mu_{2} P_{0,1}$ por lo tanto, \[\left(\lambda+\mu_{2}\right)P_{0,m}=\left(\lambda+\mu_{2}\right)\left(1-\rho_{1}\right)\rho_{2}^{m}\left(1-\rho_{2}\right).\] Entonces,

\[\mu_{2}P_{0,m+1}=\lambda\left(1-\rho_{1}\right)\rho_{2}^{m}\left(1-\rho_{2}\right)=\mu_{2}\left(1-\rho_{1}\right)\rho_{2}^{m}\left(1-\rho_{2}\right).\] De manera an\'aloga para la $\mu_1$, tenemos \[\mu_{1}P_{1,m-1}=\frac{\lambda}{\rho_{2}}\left(1-\rho_{1}\right)\rho_{2}^{m}\left(1-\rho_{2}\right).\] Luego, $\left(\lambda+\mu_{2}\right)P_{0,m}=\mu_{2}P_{0,m+1}+\mu_{1}P_{1,m-1}$, de donde, 
\[\left(\lambda+\mu_{1}+\mu_{2}\right)P_{n,m}=\left(\lambda+\mu_{1}+\mu_{2}\right)\rho^{n}\left(1-\rho_{1}\right)\rho_{2}^{m}\left(1-\rho_{2}\right),\]
por tanto
\[\mu_{2}P_{n,m+1}=\mu_{2}\rho_{2}\rho_{1}^{n}\left(1-\rho_{1}\right)\rho_{2}^{m}\left(1-\rho_{2}\right),\]

\[\mu_{1} P_{n-1,m-1}=\mu_{1}\frac{\rho_{1}}{\rho_{2}}\rho_{1}^{n}\left(1-\rho_{1}\right)\rho_{2}^{m}\left(1-\rho_{2}\right),\]

\[\lambda P_{n-1,m}=\frac{\lambda}{\rho_{1}}\rho_{1}^{n}\left(1-\rho_{1}\right)\rho_{2}^{m}\left(1-\rho_{2}\right),\]
entonces\[\left(\lambda+\mu_{1}+\mu_{2}\right)P_{n,m}=\mu_{2}P_{n,m+1}+\mu_{1} P_{n-1,m-1}+\lambda P_{n-1,m}.\]
entonces efectivamente la ecuaci\'on (\ref{Eq.8.16}) satisface las ecuaciones de balance (\ref{Eq.Balance}). El n\'umero promedio  de usuarios en el sistema, est\'a dado por
\begin{eqnarray*}
L&=&\sum_{n,m}\left(n+m\right)P_{n,m}=\sum_{n,m}nP_{n,m}+\sum_{n,m}mP_{n,m}=\sum_{n}\sum_{m}nP_{n,m}+\sum_{m}\sum_{n}mP_{n,m}\\
&=&\sum_{n}n\sum_{m}P_{n,m}+\sum_{m}m\sum_{n}P_{n,m}=\sum_{n}n\sum_{m}\rho_{1}^{n}\left(1-\rho_{1}\right)\rho_{2}^{m}\left(1-\rho_{2}\right)+\sum_{m}m\sum_{n}\rho_{1}^{n}\left(1-\rho_{1}\right)\rho_{2}^{m}\left(1-\rho_{2}\right)\\
&=&\sum_{n}n\rho_{1}^{n}\left(1-\rho_{1}\right)\sum_{m}\rho_{2}^{m}\left(1-\rho_{2}\right)+\sum_{m}m\rho_{2}^{m}\left(1-\rho_{2}\right)\sum_{n}\rho_{1}^{n}\left(1-\rho_{1}\right)=\sum_{n}n\rho_{1}^{n}\left(1-\rho_{1}\right)+\sum_{m}m\rho_{2}^{m}\left(1-\rho_{2}\right)\\
&=&\frac{\lambda}{\mu_{1}-\lambda}+\frac{\lambda}{\mu_{2}-\lambda}
\end{eqnarray*}
es decir,
\begin{equation}
L=\frac{\lambda}{\mu_{1}-\lambda}+\frac{\lambda}{\mu_{2}-\lambda}
\end{equation}

\section{Sistemas de Visitas}

Los {\emph{Sistemas de Visitas}} fueron introducidos a principios de los a\~nos 50, ver \cite{Boxma,BoonMeiWinands,LevySidi,Roubos,TakagiI,Semenova}, con un problema relacionado con las personas encargadas de la revisi\'on y reparaci\'on de m\'aquinas; m\'as adelante fueron utilizados para estudiar problemas de control de se\~nales de tr\'afico. A partir de ese momento el campo de aplicaci\'on ha crecido considerablemente, por ejemplo en: comunicaci\'on en redes de computadoras, rob\'otica, tr\'afico y transporte, manufactura, producci\'on, distribuci\'on de correo, sistema de salud p\'ublica, etc.\index{Sistemas de Visita}\\

Un modelo de colas es un modelo matem\'atico que describe la situaci\'on en la que uno o varios usuarios solicitan de un servicio a una instancia, computadora o persona. Aquellos usuarios que no son atendidos inmediatamente toman un lugar en una cola en espera de servicio. Un sistema de visitas consiste en modelos de colas conformadas por varias colas y un solo servidor que las visita en alg\'un orden para atender a los usuarios que se encuentran esperando por servicio.\index{Sistemas de Espera}

Uno de los principales objetivos de este tipo de sistemas es tratar de mejorar el desempe\~no del sistema de visitas. Una de medida de desempe\~no importante es el tiempo de respuesta del sistema, as\'i como los tiempos promedios de espera en una fila y el tiempo promedio total que tarda en ser realizada una operaci\'on completa a lo largo de todo el sistema.\index{Medidas de Desempe\~no} Algunas medidas de desempe\~no para los usuarios son los valores promedio de espera para ser atendidos, de servicio, de permanencia total en el sistema; mientras que para el servidor son los valores promedio de permanencia en una cola atendiendo, de traslado entre las colas, de duraci\'on del ciclo entre dos visitas consecutivas a la misma cola, entre otras medidas de desempe\~no estudiadas en la literatura.\\

Los sistemas de visitas pueden dividirse en dos clases:
\begin{itemize}
\item[i)] hay varios servidores y los usuarios que llegan al sistema eligen un servidor de entre los que est\'an presentes.

\item[ii)] hay uno o varios servidores que son comunes a todas las colas, estos visitan a cada una de las colas y atienden a los usuarios que est\'an presentes al momento de la visita del
servidor.
\end{itemize}

Los usuarios llegan a las colas de manera tal que los tiempos entre arribos son independientes e id\'enticamente distribuidos. En la mayor\'ia de los modelos de visitas c\'iclicas, la capacidad de almacenamiento es infinita, es decir la cola puede acomodar a una cantidad infinita de usuarios a la vez. Los tiempos de servicio en una cola son usualmente considerados como muestra de una distribuci\'on de probabilidad que caracteriza a la cola, adem\'as se acostumbra considerarlos mutuamente independientes e independientes del estado actual del sistema. \\

La ruta de atenci\'on del servidor, es el orden en el cual el servidor visita las colas determinado por un mecanismo que puede depender del estado actual del sistema (din\'amico) o puede ser independiente del estado del sistema (est\'atico). El mecanismo m\'as utilizado es el c\'iclico. Para modelar sistemas en los cuales ciertas colas son visitadas con mayor frecuencia que otras, las colas c\'iclicas se han extendido a colas peri\'odicas, en las cuales el servidor visita la cola conforme a una orden de servicio de longitud finita. \index{Colas C\'iclicas}\\

El {\em orden de visita} se entiende como la regla utilizada por el servidor para elegir la pr\'oxima cola. Este servicio puede ser din\'amico o est\'atico:\index{Orden de visita}

\begin{itemize}
\item[i)] Para el caso {\em est\'atico} la regla permanece invariante a lo largo del curso de la operaci\'on del sistema.

\item[ii)] Para el caso {\em din\'amico} la cola que se elige para servicio en el momento depende de un conocimiento total o parcial del estado del sistema.
\end{itemize}

Dentro de los ordenes de tipo est\'atico hay varios, los m\'as comunes son:

\begin{itemize}
\item[i)] {\em c\'iclico}: Si denotamos por $\left\{Q_{i}\right\}_{i=1}^{N}$ al conjunto de colas a las cuales el servidor visita en el orden \[Q_{1},Q_{2},\ldots,Q_{N},Q_{1},Q_{2},\ldots,Q_{N}.\]

\item[ii)] {\em peri\'odico}: el servidor visita las colas en el orden:
\[Q_{T\left(1\right)},Q_{T\left(2\right)},\ldots,Q_{T\left(M\right)},Q_{T\left(1\right)},\ldots,Q_{T\left(M\right)}\]
caracterizada por una tabla de visitas
\[\left(T\left(1\right),T\left(2\right),\ldots,T\left(M\right)\right),\]
con $M\geq N$, $T\left(i\right)\in\left\{1,2,\ldots,N\right\}$ e $i=\overline{1,M}$. Hay un caso especial, {\em colas tipo elevador} donde las colas son atendidas en el orden \[Q_{1},Q_{2},\ldots,Q_{N},Q_{1},Q_{2},\ldots,Q_{N-1},Q_{N},Q_{N-1},\ldots,Q_{1},\].

\item[iii)] {\em aleatorio}: la cola $Q_{i}$ es elegida para ser atendida con probabilidad $p_{i}$, $i=\overline{1,N}$, $\sum_{i=1}^{N}p_{i}=1$. Una posible variaci\'on es que despu\'es de atender $Q_{i}$ el servidor se desplaza a $Q_{j}$ con probabilidad $p_{ij}$, con $i,j=\overline{1,N}$, $\sum_{j=1}^{N}p_{ij}=1$, para $i=\overline{1,N}$.
\end{itemize}

El servidor usualmente incurrir\'a en tiempos de traslado para ir de una cola a otra. Un sistema de visitas puede expresarse en un par de par\'ametros: el n\'umero de colas, que usualmente se denotar\'a por $N$, y el tr\'afico caracter\'istico de las colas, que consiste de los procesos de arribo y los procesos de servicio caracteriza a estos sistemas.\\

La disciplina de servicio especifica el n\'umero de usuarios que son atendidos durante la visita del servidor a la cola; estas pueden ser clasificadas en l\'imite de usuarios atendidos y en usuarios atendidos en un tiempo l\'imite, poniendo restricciones en la cantidad de tiempo utilizado por el servidor en una visita a la cola. Alternativamente pueden ser clasificadas en pol\'iticas exhaustivas y pol\'iticas cerradas, dependiendo en si los usuarios que llegaron a la cola mientras el servidor estaba dando servicio son candidatos para ser atendidos por el servidor que se encuentra en la cola dando servicio. En la pol\'itica exhaustiva estos usuarios son candidatos para ser atendidos mientras que en la cerrada no lo son. De estas dos pol\'iticas se han creado h\'ibridos los cuales pueden revisarse en \cite{BoonMeiWinands}.\index{Disciplina de Servicio}\\

La disciplina de la cola especifica el orden en el cual los usuarios presentes en la cola son atendidos. La m\'as com\'un es la {\em First-In-First-Served}. Las pol\'iticas m\'as comunes son las de tipo exhaustivo que consiste en que el servidor continuar\'a trabajando hasta que la cola quede vac\'ia; y la pol\'itica cerrada, bajo la cual ser\'an atendidos exactamente aquellos que estaban presentes al momento en que lleg\'o el servidor a la cola. \index{Pol\'itica de Servicio}\\

Las pol\'iticas de servicio deben de satisfacer las siguientes propiedades:
\begin{itemize}
\item[i)] No dependen de los procesos de servicio anteriores.
\item[ii)] La selecci\'on de los usuarios para ser atendidos es independiente del tiempo de servicio requerido  y de los posibles arribos futuros.
\item[iii)] las pol{\'\i}ticas de servicio que son aplicadas, es decir, el n\'umero de usuarios en la cola que ser{\'a}n atendidos durante la visita del servidor a la misma; \'estas pueden ser clasificadas por la cantidad de usuarios atendidos y por el n\'umero de usuarios atendidos en un intervalo de tiempo determinado. Las principales pol\'iticas de servicio para las cuales se han desarrollado aplicaciones son: la exhaustiva, la cerrada y la $k$-l\'imite, ver \cite{LevySidi, TakagiI, Semenova}. De estas pol\'iticas se han creado h\'ibridos los cuales pueden revisarse en Boon and Van der Mei \cite{BoonMeiWinands}.

\item[iv)] Una pol{\'\i}tica de servicio es asignada a cada etapa independiente de la cola que se est{\'a} atendiendo, no necesariamente es la misma para todas las etapas.
\item[v)] El servidor da servicio de manera constante.

\item[vi)] La pol\'itica de servicio se asume mon\'otona (ver
\cite{Stability}).

\end{itemize}

Las principales pol\'iticas deterministas de servicio son:\index{Pol\'iticas Deterministas}
\begin{itemize}

\item[i)] {\em Cerrada} donde solamente los usuarios presentes al comienzo de la etapa son considerados para ser atendidos.

\item[ii)] {\em Exhaustiva} en la que tanto los usuarios presentes al comienzo de la etapa como los que arriban   mientras se est\'a dando servicio son considerados para ser atendidos.

\item[iii)] $k_{i}$-limited: el n\'umero de usuarios por atender en la cola $i$ est\' acotado por $k_{i}$.

\item[iv)] {\em tiempo limitado} la cola es atendida solo por un periodo de tiempo fijo.
\end{itemize}

\begin{Note}
\begin{itemize}
\item[a) ] Una etapa es el periodo de tiempo durante el cual el
servidor atiende de manera continua en una sola cola.

\item[b) ] Un ciclo  es el periodo necesario para terminar $l$ etapas.
\end{itemize}
\end{Note}

Boxma y Groenendijk \cite{Boxma2} enuncian la Ley de Pseudo-Conservaci\'on para la pol\'itica exhaustiva como\index{Ley de Pseudo-Conservaci\'on}

\begin{equation}\label{LPCPE}
\sum_{i=1}^{N}\rho_{i}\esp
W_{i}=\rho\frac{\sum_{i=1}^{N}\lambda_{i}\esp\left[\delta_{i}^{(2)}\left(1\right)\right]}{2\left(1-\rho\right)}+\rho\frac{\delta^{(2)}}{2\delta}+\frac{\delta}{2\left(1-\rho\right)}\left[\rho^{2}-\sum_{i=1}^{N}\rho_{i}^{2}\right],
\end{equation}

donde $\delta=\sum_{i=1}^{N}\delta_{i}\left(1\right)$ y
$\delta_{i}^{(2)}$ denota el segundo momento de los tiempos de traslado entre colas del servidor, $\delta^{(2)}$ es el segundo momento de los tiempos de traslado entre las colas de todo el sistema, finalmente $\rho=\sum_{i=1}^{N}\rho_{i}$. Por otro lado, se tiene que

\begin{equation}\label{Eq.Tiempo.Espera}\index{Tiempos de Espera}
\esp W_{i}=\frac{\esp I_{i}^{2}}{2\esp
I_{i}}+\frac{\lambda_{i}\esp\left[\eta_{i}^{(2)}\left(1\right)\right]}{2\left(1-\rho_{i}\right)},
\end{equation}

con $I_{i}$ definido como el peri\'odo de intervisita, es decir el tiempo entre una salida y el pr\'oximo arribo del servidor a la cola $Q_{i}$, dado por $I_{i}=C_{i}-V_{i}$, donde $C_{i}$ es la longitud del ciclo, definido como el tiempo entre dos instantes de visita consecutivos a la cola $Q_{i}$ y $V_{i}$ es el periodo de visita, definido como el tiempo que el servidor utiliza en atender a los usuarios de la cola $Q_{i}$.\index{Periodo de Intervisita}
\begin{equation}\label{Eq.Periodo.Intervisita}
\esp
I_{i}=\frac{\left(1-\rho_{i}\right)}{1-\rho}\sum_{i=1}^{N}\esp\left[\delta_{i}\left(1\right)\right],
\end{equation}

con

\begin{equation}\label{SdoMomento.Periodo.Intervisita}
\esp
I_{i}^{2}=\esp\left[\delta_{i-1}^{(2)}\left(1\right)\right]-\left(\esp\left[\delta_{i-1}\left(1\right)\right]\right)^{2}+
\frac{1-\rho_{i}}{\rho_{i}}\sum_{j=1,j\neq i}^{N}r_{ij}+\left(\esp
I_{i}\right)^{2},
\end{equation}

donde el conjunto de valores $\left\{r_{ij}:i,j=1,2,\ldots,N\right\}$ representan la covarianza del tiempo para las colas $i$ y $j$; para sistemas con servicio exhaustivo, el tiempo de estaci\'on para la cola $i$ se define como el intervalo de tiempo entre instantes sucesivos cuando el servidor abandona la cola $i-1$ y la cola $i$. Hideaki Takagi \cite{Takagi} proporciona expresiones cerradas para calcular $r_{ij}$, \'estas implican resolver un sistema de $N^{2}$ ecuaciones lineales;

\begin{eqnarray}\label{Eq.Cov.TT}
\begin{array}{l}
r_{ij}=\frac{\rho_{i}}{1-\rho_{i}}\left(\sum_{m=i+1}^{N}r_{jm}+\sum_{m=1}^{j-1}r_{jm}+\sum_{m=j}^{i-1}r_{jm}\right),\textrm{
}j<i,\\
r_{ij}=\frac{\rho_{i}}{1-\rho_{i}}\left(\sum_{m=i+1}^{j-1}r_{jm}+\sum_{m=j}^{N}r_{jm}+\sum_{m=1}^{i-1}r_{jm}\right),\textrm{
}j>i,\\
r_{ij}=\frac{\esp\left[\delta_{i-1}^{(2)}\left(1\right)\right]-\left(\esp\left[\delta_{i-1}\left(1\right)\right]\right)^{2}}
{\left(1-\rho_{i}\right)^{2}}+\frac{\lambda_{i}\esp\left[\eta_{i}\left(1\right)^{(2)}\right]}{\left(1-\rho_{i}\right)^{3}}+\frac{\rho_{i}}{1-\rho_{i}}\sum_{j=i,j=1}^{N}r_{ij}.
\end{array}
\end{eqnarray}\index{Pol\'itica de Servicio Exhaustiva}

Para el caso de la Pol\'itica Cerrada la Ley de Pseudo-Conservaci\'on se expresa en los siguientes t\'erminos.
\begin{equation}\label{LPCPG}
\sum_{i=1}^{N}\rho_{i}\esp
W_{i}=\rho\frac{\sum_{i=1}^{N}\lambda_{i}\esp\left[\delta_{i}\left(1\right)^{(2)}\right]}{2\left(1-\rho\right)}+\rho\frac{\delta^{(2)}}{2\delta}+\frac{\delta}{2\left(1-\rho\right)}\left[\rho^{2}+\sum_{i=1}^{N}\rho_{i}^{2}\right],
\end{equation}
el tiempo de espera promedio para los usuarios en la cola $Q_{1}$ se puede determinar por medio de
\begin{equation}\label{Eq.Tiempo.Espera.Gated}
\esp W_{i}=\frac{\left(1+\rho_{i}\right)\esp C_{i}^{2}}{2\esp
C_{i}},
\end{equation}\index{Tiempos de Espera}
donde $C_{i}$ denota la longitud del ciclo para la cola $Q_{i}$, definida como el tiempo entre dos instantes consecutivos de visita en $Q_{i}$, cuyo segundo momento est\'a dado por

\begin{equation}\label{Eq.Periodo.Intervisita.Gated}
\esp C_{i}^{2}=\frac{1}{\rho_{i}}\sum_{j=1,j\neq
i}^{N}r_{ij}+\sum_{j=1}^{N}r_{ij}+\left(\esp C\right)^{2},
\end{equation}\index{Periodo de Intervisita}
con
\begin{eqnarray*}
\esp C=\frac{\delta}{1-\rho},
\end{eqnarray*}
donde $r_{ij}$ representa la covarianza del tiempo de estaci\'on para las colas $i$ y $j$, pero el tiempo de estaci\'on para la cola $i$ para la pol\'itica cerrada se define como el intervalo de tiempo entre instantes sucesivos cuando el servidor visita la cola $i$ y la cola $i+1$. El conjunto $\left\{r_{ij}:i,j=1,2,\ldots,N\right\}$ se calcula resolviendo un
sistema de $N^{2}$ ecuaciones lineales

\begin{eqnarray}\label{Eq.Cov.TT.Gated}
\begin{array}{l}
r_{ij}=\rho_{i}\left(\sum_{m=i}^{N}r_{jm}+\sum_{m=1}^{j-1}r_{jm}+\sum_{m=j}^{i-1}r_{mj}\right),\textrm{
}j<i,\\
r_{ij}=\rho_{i}\left(\sum_{m=i}^{j-1}r_{jm}+\sum_{m=j}^{N}r_{jm}+\sum_{m=1}^{i-1}r_{mj}\right),\textrm{
}j>i,\\
r_{ij}=r_{i-1}^{(2)}-\left(r_{i-1}^{(1)}\right)^{2}+\lambda_{i}b_{i}^{(2)}\esp
C_{i}+\rho_{i}\sum_{j=1,j\neq
i}^{N}r_{ij}+\rho_{i}^{2}\sum_{i=j,j=1}^{N}r_{ij}.
\end{array}
\end{eqnarray}

Finalmente, Takagi \cite{TakagiI} proponen una aproximaci\'on para los tiempos de espera de los usuarios en cada una de las colas:
\begin{eqnarray}
\begin{array}{l}
\sum_{i=1}^{N}\frac{\rho_{i}}{\rho}\left(1-\frac{\lambda_{i}\delta}{1-\rho}\right)\esp\left[W_{i}\right]=\sum_{i=1}^{N}\frac{\lambda_{i}\esp\left[\eta_{i}\left(1\right)^{(2)}\right]}{2\left(1-\rho\right)}\\
+\frac{\sum_{i=1}^{N}\esp\left[\delta_{i}^{2}\right]-\left(\esp\left[\delta_{i}\left(1\right)\right]\right)^{2}}{2\delta}+\frac{\delta\left(\rho-\sum_{i=1}^{N}\rho_{i}^{2}\right)}{2\rho\left(1-\rho\right)}+\frac{\delta\sum_{i=1}^{N}\rho_{i}^{2}}{\rho\left(1-\rho\right)},
\end{array}
\end{eqnarray}
entonces
\begin{eqnarray}\label{LPCPKL}
\begin{array}{l}
\esp
W_{i}\cong\frac{1-\rho+\rho_{i}}{1-\rho-\lambda_{i}\delta}\times\frac{1-\rho}{\rho\left(1-\rho\right)+\sum_{i=1}^{N}\rho_{i}^{2}}\\
\times\left[\frac{\rho}{2\left(1-\rho\right)}\sum_{i=1}^{N}\lambda_{i}\esp\left[\eta_{i}\left(1\right)^{(2)}\right]+\frac{\rho\Delta^{2}}{2\delta}+\frac{\delta}{2\left(1-\rho\right)}\sum_{i=1}^{N}\rho_{i}\left(1+\rho_{i}\right).\right]
\end{array}
\end{eqnarray}
donde $\Delta^{2}=\sum_{i=1}^{N}\delta_{i}^{2}$.

 \section{Funci\'on Generadora de Probabilidades}

\begin{Teo}[Teorema de Continuidad]\index{Teorema de Continuidad}
Sup\'ongase que $\left\{X_{n},n=1,2,3,\ldots\right\}$ son variables aleatorias finitas, no negativas con valores enteros tales que $P\left(X_{n}=k\right)=p_{k}^{(n)}$, para $n=1,2,3,\ldots$, $k=0,1,2,\ldots$, con $\sum_{k=0}^{\infty}p_{k}^{(n)}=1$, para $n=1,2,3,\ldots$. Sea $g_{n}$ la Funci\'on Generadora de Probabilidades (FGP)  para la variable aleatoria $X_{n}$. Entonces existe una sucesi\'on $\left\{p_{k}\right\}$ tal que 
\begin{eqnarray}
lim_{n\rightarrow\infty}p_{k}^{(n)}=p_{k}\textrm{ para }0<s<1.
\end{eqnarray}

En este caso, $g\left(s\right)=\sum_{k=0}^{\infty}s^{k}p_{k}$. Adem\'as
\begin{eqnarray}
\sum_{k=0}^{\infty}p_{k}=1\textrm{ si y s\'olo si
}lim_{s\uparrow1}g\left(s\right)=1.
\end{eqnarray}
\end{Teo}

\begin{Teo}
Sea $N$ una variable aleatoria con valores enteros no negativos finita tal que $P\left(N=k\right)=p_{k}$, para $k=0,1,2,\ldots$, y 
\begin{eqnarray}
\sum_{k=0}^{\infty}p_{k}=P\left(N<\infty\right)=1.
\end{eqnarray} 

Sea $\Phi$ la FGP de $N$ tal que
\begin{eqnarray}
g\left(s\right)=\esp\left[s^{N}\right]=\sum_{k=0}^{\infty}s^{k}p_{k},
\end{eqnarray}
 con $g\left(1\right)=1$. Si $0\leq p_{1}\leq1$ y 
 \begin{eqnarray}
\esp\left[N\right]=g^{'}\left(1\right)\leq1,
\end{eqnarray}
 entonces no existe soluci\'on  de la ecuaci\'on $g\left(s\right)=s$ en el intervalo $\left[0,1\right)$. Si $\esp\left[N\right]=g^{'}\left(1\right)>1$, lo cual implica que $0\leq p_{1}<1$, entonces existe una \'unica soluci\'on de la ecuaci\'on $g\left(s\right)=s$ en el intervalo $\left[0,1\right)$.
\end{Teo}

\begin{Teo}
Si $X$ y $Y$ tienen PGF $G_{X}$ y $G_{Y}$ respectivamente, entonces,\[G_{X}\left(s\right)=G_{Y}\left(s\right)\] para toda $s$, s\'i y s\'olo s\'i
\begin{eqnarray}
P\left(X=k\right))=P\left(Y=k\right),
\end{eqnarray}
para toda $k=0,1,\ldots,$., es decir, si y s\'olo si $X$ y $Y$ tienen la misma distribuci\'on de probabilidad.
\end{Teo}

\begin{Teo}
Para cada $n$ fijo, sea la sucesi\'on de probabilidades $\left\{a_{0,n},a_{1,n},\ldots,\right\}$, tales que $a_{k,n}\geq0$ para toda $k=0,1,2,\ldots,$ y $\sum_{k\geq0}a_{k,n}=1$, y sea $G_{n}\left(s\right)$ la correspondiente funci\'on generadora, $G_{n}\left(s\right)=\sum_{k\geq0}a_{k,n}s^{k}$. De modo que para cada valor fijo de $k$
\begin{eqnarray}
lim_{n\rightarrow\infty}a_{k,n}=a_{k},
\end{eqnarray}
es decir converge en distribuci\'on, es necesario y suficiente que para cada valor fijo $s\in\left[0,\right)$,

\begin{eqnarray}
lim_{n\rightarrow\infty}G_{n}\left(s\right)=G\left(s\right),
\end{eqnarray}
donde $G\left(s\right)=\sum_{k\geq0}p_{k}s^{k}$, para cualquier la funci\'on generadora del l\'imite de la sucesi\'on.
\end{Teo}

\begin{Teo}[Teorema de Abel]\index{Teorema de Abel}
Sea $G\left(s\right)=\sum_{k\geq0}a_{k}s^{k}$ para cualquier $\left\{p_{0},p_{1},\ldots,\right\}$, tales que $p_{k}\geq0$ para toda $k=0,1,2,\ldots,$. Entonces $G\left(s\right)$ es continua por la derecha en $s=1$, es decir
\begin{eqnarray}
lim_{s\uparrow1}G\left(s\right)=\sum_{k\geq0}p_{k}=G\left(s\right),
\end{eqnarray}
sin importar si la suma es finita o no.
\end{Teo}

\begin{Note}
El radio de Convergencia para cualquier FGP es $R\geq1$, entonces, el Teorema de Abel nos dice que a\'un en el peor escenario, cuando $R=1$, a\'un se puede confiar en que la FGP ser\'a continua en $s=1$, en contraste, no se puede asegurar que la FGP ser\'a continua en el l\'imite inferior $-R$, puesto que la FGP es sim\'etrica alrededor del cero: la FGP converge para todo $s\in\left(-R,R\right)$, y no lo hace para $s<-R$ o $s>R$. Adem\'as nos dice que podemos escribir $G_{X}\left(1\right)$ como una abreviaci\'on de $lim_{s\uparrow1}G_{X}\left(s\right)$.
\end{Note}

Entonces si suponemos que la diferenciaci\'on t\'ermino a t\'ermino est\'a permitida, entonces

\begin{eqnarray}
G_{X}^{'}\left(s\right)&=&\sum_{x=1}^{\infty}xs^{x-1}p_{x}
\end{eqnarray}

el Teorema de Abel nos dice que
\begin{eqnarray}
\begin{array}{l}
\esp\left(X\right]=\lim_{s\uparrow1}G_{X}^{'}\left(s\right):\\
\esp\left[X\right]=\sum_{x=1}^{\infty}xp_{x}=G_{X}^{'}\left(1\right)=\lim_{s\uparrow1}G_{X}^{'}\left(s\right),
\end{array}
\end{eqnarray}
dado que el Teorema de Abel se aplica a
\begin{eqnarray}
G_{X}^{'}\left(s\right)&=&\sum_{x=1}^{\infty}xs^{x-1}p_{x},
\end{eqnarray}
estableciendo as\'i que $G_{X}^{'}\left(s\right)$ es continua en $s=1$. Sin el Teorema de Abel no se podr\'ia asegurar que el l\'imite de $G_{X}^{'}\left(s\right)$ conforme $s\uparrow1$ sea la respuesta correcta para $\esp\left[X\right]$.

\begin{Note}
La FGP converge para todo $|s|<R$, para alg\'un $R$. De hecho la FGP converge absolutamente si $|s|<R$. La FGP adem\'as converge uniformemente en conjuntos de la forma $\left\{s:|s|<R^{'}\right\}$, donde $R^{'}<R$, es decir, $\forall\epsilon>0, \exists n_{0}\in\ent$ tal que $\forall s$, con $|s|<R^{'}$, y $\forall n\geq n_{0}$,
\begin{eqnarray}
|\sum_{x=0}^{n}s^{x}\prob\left(X=x\right)-G_{X}\left(s\right)|<\epsilon.
\end{eqnarray}
De hecho, la convergencia uniforme es la que nos permite diferenciar t\'ermino a t\'ermino:
\begin{eqnarray}
G_{X}\left(s\right)=\esp\left[s^{X}\right]=\sum_{x=0}^{\infty}s^{x}\prob\left(X=x\right),
\end{eqnarray}
y sea $s<R$. Entonces se tiene lo siguiente: 
\begin{eqnarray}
G_{X}^{'}\left(s\right)=\frac{d}{ds}\left(\sum_{x=0}^{\infty}s^{x}\prob\left(X=x\right)\right)=\sum_{x=0}^{\infty}\frac{d}{ds}\left(s^{x}\prob\left(X=x\right)\right)=\sum_{x=0}^{n}xs^{x-1}\prob\left(X=x\right).
\end{eqnarray}

\begin{eqnarray}
\int_{a}^{b}G_{X}\left(s\right)ds&=&\int_{a}^{b}\left(\sum_{x=0}^{\infty}s^{x}\prob\left(X=x\right)\right)ds=\sum_{x=0}^{\infty}\left(\int_{a}^{b}s^{x}\prob\left(X=x\right)ds\right)=\sum_{x=0}^{\infty}\frac{s^{x+1}}{x+1}\prob\left(X=x\right),
\end{eqnarray}
para $-R<a<b<R$.

\end{Note}

\begin{Teo}[Teorema de Convergencia Mon\'otona para FGP]\index{Teorema de Convergencia Mon\'otona} Sean $X$ y $X_{n}$ variables aleatorias no negativas, con valores en los enteros, finitas, tales que
\begin{eqnarray*}
lim_{n\rightarrow\infty}G_{X_{n}}\left(s\right)&=&G_{X}\left(s\right)\textrm{, para }0\leq s\leq1,
\end{eqnarray*}
entonces
\begin{eqnarray*}
lim_{n\rightarrow\infty}P\left(X_{n}=k\right)=P\left(X=k\right)\textrm{, para }k=0,1,2,\ldots.
\end{eqnarray*}

\end{Teo}

El teorema anterior requiere del siguiente lema:
\begin{Lema}
Sean $a_{n,k}\in\ent^{+}$, $n\in\nat$ constantes no negativas con $\sum_{k\geq0}a_{k,n}\leq1$. Sup\'ongase que para $0\leq s\leq1$, se tiene
\begin{eqnarray}
a_{n}\left(s\right)&=&\sum_{k=0}^{\infty}a_{k,n}s^{k}\rightarrow
a\left(s\right)=\sum_{k=0}^{\infty}a_{k}s^{k}.
\end{eqnarray}
Entonces
\begin{eqnarray}
a_{0,n}\rightarrow a_{0}.
\end{eqnarray}
\end{Lema}

Consideremos un sistema que consta de \'unicamente un servidor y una sola cola, a la cual los usuarios arriban conforme a un proceso Poisson cuya tasa promedio de llegada es $1/\lambda$; la tasa promedio con la cual el servidor da servicio es $1/\mu$, adem\'as, los tiempos entre arribos y los tiempos de servicio son independientes entre s\'i. Se define la carga de tr\'afico $\rho:=\frac{\lambda}{\mu}$, para este modelo existe un teorema que nos dice la relaci\'on que hay entre el valor de $\rho$ y la estabilidad de la cola:\index{Proceso Poisson}

\begin{Prop}\index{Cola $M/M/1$}\index{Cola $M/M/1$}
La cola $M/M/1$ con carga de tr\'afico $\rho$, es estable si y s\'olo si $\rho<1$.
\end{Prop}

Este teorema nos permite determinar las principales medidas de desempe\~no: Tiempo de espera en el sistema, $W$, el n\'umero esperado de clientes en el sistema, $L$, adem\'as de los tiempos promedio e espera tanto en la cola como de servicio, $s$ representa el tiempo de servicio para un cliente:

\begin{eqnarray}
\begin{array}{ll}
 L=\frac{\rho}{1-\rho},& W=\frac{1}{\mu-\lambda},\\
W_{q}=\esp\left[s\right]\frac{\rho}{1-\rho}\textrm{,  y }&  L_{q}=\frac{\rho^{2}}{1-\rho}.
\end{array}
\end{eqnarray}

Esta es la idea general, poder determinar la principales medidas de desempe\~no para un sistema de colas o sistema de visitas, para este fin es necesario realizar los siguientes supuestos. En teor\'ia de colas hay casos particulares, para los cuales es posible determinar espec\'ificamente medidas de desempe\~no del sistema bajo condiciones de estabilidad, tales como los tiempos promedio de espera y de servicio, tanto en el sistema como en cada una de las colas. Se considerar\'an intervalos de tiempo de la forma $\left[t,t+1\right]$, adem\'as sup\'ongase que los usuarios arriban por paquetes de manera independiente del resto de las colas. Se define el grupo de usuarios que llegan a cada una de las colas del sistema 1, caracterizadas por $Q_{1}$ y $Q_{2}$ respectivamente, en el intervalo de tiempo $\left[t,t+1\right]$ por $X_{1}\left(t\right),X_{2}\left(t\right)$.\\

Para cada uno de los procesos anteriores se define su Funci\'on Generadora de Probabilidades (FGP):

\begin{eqnarray}
\begin{array}{cc}
P_{1}\left(z_{1}\right)=\esp\left[z_{1}^{X_{1}\left(t\right)}\right], & P_{2}\left(z_{2}\right)=\esp\left[z_{2}^{X_{2}\left(t\right)}\right].
\end{array}
\end{eqnarray}

Con primer momento definidos por

\begin{eqnarray}
\begin{array}{ll}
\mu_{1}=\esp\left[X_{1}\left(t\right)\right]=P_{1}^{(1)}\left(1\right),& \mu_{2}=\esp\left[X_{2}\left(t\right)\right]=P_{2}^{(1)}\left(1\right).
\end{array}
\end{eqnarray}

En lo que respecta al servidor, en t\'erminos de los tiempos de visita a cada una de las colas, se denotar\'an por $\tau_{1},\tau_{2}$ para $Q_{1},Q_{2}$ respectivamente; y a los tiempos en que el servidor termina de atender en las colas $Q_{1},Q_{2}$, se les denotar\'a por $\overline{\tau}_{1},\overline{\tau}_{2}$ respectivamente. Entonces, los tiempos de servicio est\'an dados por las diferencias $\overline{\tau}_{1}-\tau_{1},\overline{\tau}_{2}-\tau_{2}$ para $Q_{1},Q_{2}$. \\

An\'alogamente los tiempos de traslado del servidor desde el momento en que termina de atender a una cola y llega a la siguiente para comenzar a dar servicio est\'an dados por $\tau_{2}-\overline{\tau}_{1},\tau_{1}-\overline{\tau}_{2}$. La FGP para estos tiempos de traslado est\'an dados por

\begin{eqnarray}
\begin{array}{cc}
R_{1}\left(z_{1}\right)=\esp\left[z_{1}^{\tau_{2}-\overline{\tau}_{1}}\right], & R_{2}\left(z_{2}\right)=\esp\left[z_{2}^{\tau_{1}-\overline{\tau}_{2}}\right],
\end{array}
\end{eqnarray}

y al igual que como se hizo con anterioridad

\begin{eqnarray}
\begin{array}{cc}
r_{1}=R_{1}^{(1)}\left(1\right)=\esp\left[\tau_{2}-\overline{\tau}_{1}\right], & r_{2}=R_{2}^{(1)}\left(1\right)=\esp\left[\tau_{1}-\overline{\tau}_{2}\right].
\end{array}
\end{eqnarray}
Sean $\alpha_{1},\alpha_{2}$ el n\'umero de usuarios que arriban
en grupo a la cola $Q_{1}$ y $Q_{2}$ respectivamente. Sus FGP's
est\'an definidas como

\begin{eqnarray}
\begin{array}{cc}
A_{1}\left(z\right)=\esp\left[z^{\alpha_{1}\left(t\right)}\right],&
A_{2}\left(z\right)=\esp\left[z^{\alpha_{2}\left(t\right)}\right].
\end{array}
\end{eqnarray}

Su primer momento est\'a dado por

\begin{eqnarray}
\begin{array}{cc}
\lambda_{1}=\esp\left[\alpha_{1}\left(t\right)\right]=A_{1}^{(1)}\left(1\right),&
\lambda_{2}=\esp\left[\alpha_{2}\left(t\right)\right]=A_{2}^{(1)}\left(1\right).
\end{array}
\end{eqnarray}

Sean $\beta_{1},\beta_{2}$ el n\'umero de usuarios que arriban en el grupo $\alpha_{1},\alpha_{2}$ a la cola $Q_{1}$ y $Q_{2}$, respectivamente, de igual manera se definen sus FGP's

\begin{eqnarray}
\begin{array}{cc}
B_{1}\left(z\right)=\esp\left[z^{\beta_{1}\left(t\right)}\right],&
B_{2}\left(z\right)=\esp\left[z^{\beta_{2}\left(t\right)}\right],
\end{array}
\end{eqnarray}

con

\begin{eqnarray}
\begin{array}{cc}
b_{1}=\esp\left[\beta_{1}\left(t\right)\right]=B_{1}^{(1)}\left(1\right),&
b_{2}=\esp\left[\beta_{2}\left(t\right)\right]=B_{2}^{(1)}\left(1\right).
\end{array}
\end{eqnarray}

La distribuci\'on para el n\'umero de grupos que arriban al sistema en cada una de las colas se definen por:

\begin{eqnarray}
\begin{array}{cc}
P_{1}\left(z_{1}\right)=A_{1}\left[B_{1}\left(z_{1}\right)\right]=\esp\left[B_{1}\left(z_{1}\right)^{\alpha_{1}\left(t\right)}\right],& P_{2}\left(z_{1}\right)=A_{1}\left[B_{1}\left(z_{1}\right)\right]=\esp\left[B_{1}\left(z_{1}\right)^{\alpha_{1}\left(t\right)}\right],
\end{array}
\end{eqnarray}

entonces

\begin{eqnarray}
\begin{array}{l}
P_{1}^{(1)}\left(1\right)=\esp\left[\alpha_{1}\left(t\right)B_{1}^{(1)}\left(1\right)\right]=B_{1}^{(1)}\left(1\right)\esp\left[\alpha_{1}\left(t\right)\right]=\lambda_{1}b_{1}\\
P_{2}^{(1)}\left(1\right)=\esp\left[\alpha_{2}\left(t\right)B_{2}^{(1)}\left(1\right)\right]=B_{2}^{(1)}\left(1\right)\esp\left[\alpha_{2}\left(t\right)\right]=\lambda_{2}b_{2}.
\end{array}
\end{eqnarray}

De lo desarrollado hasta ahora se tiene lo siguiente

\begin{eqnarray*}
\begin{array}{l}
\esp\left[z_{1}^{L_{1}\left(\overline{\tau}_{1}\right)}z_{2}^{L_{2}\left(\overline{\tau}_{1}\right)}\right]=\esp\left[z_{2}^{L_{2}\left(\overline{\tau}_{1}\right)}\right]=\esp\left[z_{2}^{L_{2}\left(\tau_{1}\right)+X_{2}\left(\overline{\tau}_{1}-\tau_{1}\right)}\right]=\esp\left[\left\{z_{2}^{L_{2}\left(\tau_{1}\right)}\right\}\left\{z_{2}^{X_{2}\left(\overline{\tau}_{1}-\tau_{1}\right)}\right\}\right]\\
=\esp\left[\left\{z_{2}^{L_{2}\left(\tau_{1}\right)}\right\}\left\{P_{2}\left(z_{2}\right)\right\}^{\overline{\tau}_{1}-\tau_{1}}\right]=\esp\left[\left\{z_{2}^{L_{2}\left(\tau_{1}\right)}\right\}\left\{\theta_{1}\left(P_{2}\left(z_{2}\right)\right)\right\}^{L_{1}\left(\tau_{1}\right)}\right]=F_{1}\left(\theta_{1}\left(P_{2}\left(z_{2}\right)\right),z_{2}\right),
\end{array}
\end{eqnarray*}

es decir 
\begin{equation}\label{Eq.base.F1}
\esp\left[z_{1}^{L_{1}\left(\overline{\tau}_{1}\right)}z_{2}^{L_{2}\left(\overline{\tau}_{1}\right)}\right]=F_{1}\left(\theta_{1}\left(P_{2}\left(z_{2}\right)\right),z_{2}\right).
\end{equation}

Procediendo de manera an\'aloga para $\overline{\tau}_{2}$:

\begin{eqnarray*}
\begin{array}{l}
\esp\left[z_{1}^{L_{1}\left(\overline{\tau}_{2}\right)}z_{2}^{L_{2}\left(\overline{\tau}_{2}\right)}\right]=\esp\left[z_{1}^{L_{1}\left(\overline{\tau}_{2}\right)}\right]=\esp\left[z_{1}^{L_{1}\left(\tau_{2}\right)+X_{1}\left(\overline{\tau}_{2}-\tau_{2}\right)}\right]=\esp\left[\left\{z_{1}^{L_{1}\left(\tau_{2}\right)}\right\}\left\{z_{1}^{X_{1}\left(\overline{\tau}_{2}-\tau_{2}\right)}\right\}\right]\\
=\esp\left[\left\{z_{1}^{L_{1}\left(\tau_{2}\right)}\right\}\left\{P_{1}\left(z_{1}\right)\right\}^{\overline{\tau}_{2}-\tau_{2}}\right]=\esp\left[\left\{z_{1}^{L_{1}\left(\tau_{2}\right)}\right\}\left\{\theta_{2}\left(P_{1}\left(z_{1}\right)\right)\right\}^{L_{2}\left(\tau_{2}\right)}\right]=F_{2}\left(z_{1},\theta_{2}\left(P_{1}\left(z_{1}\right)\right)\right),
\end{array}
\end{eqnarray*}
por tanto
\begin{equation}\label{Eq.PGF.Conjunta.Tau2}
\esp\left[z_{1}^{L_{1}\left(\overline{\tau}_{2}\right)}z_{2}^{L_{2}\left(\overline{\tau}_{2}\right)}\right]=F_{2}\left(z_{1},\theta_{2}\left(P_{1}\left(z_{1}\right)\right)\right).
\end{equation}

Ahora, para el intervalo de tiempo
$\left[\overline{\tau}_{1},\tau_{2}\right]$ y $\left[\overline{\tau}_{2},\tau_{1}\right]$, los arribos de los usuarios modifican el n\'umero de usuarios que llegan a las colas, es decir, los procesos
$L_{1}\left(t\right)$ y $L_{2}\left(t\right)$. La FGP para el n\'umero de arribos a todas las estaciones durante el intervalo $\left[\overline{\tau}_{1},\tau_{2}\right]$  cuya distribuci\'on est\'a especificada por la distribuci\'on compuesta $R_{1}\left(\mathbf{z}\right),R_{2}\left(\mathbf{z}\right)$:

\begin{eqnarray*}
R_{1}\left(\mathbf{z}\right)=R_{1}\left(\prod_{i=1}^{2}P\left(z_{i}\right)\right)=\esp\left[\left\{\prod_{i=1}^{2}P\left(z_{i}\right)\right\}^{\tau_{2}-\overline{\tau}_{1}}\right]\\
R_{2}\left(\mathbf{z}\right)=R_{2}\left(\prod_{i=1}^{2}P\left(z_{i}\right)\right)=\esp\left[\left\{\prod_{i=1}^{2}P\left(z_{i}\right)\right\}^{\tau_{1}-\overline{\tau}_{2}}\right].
\end{eqnarray*}

Dado que los eventos en
$\left[\tau_{1},\overline{\tau}_{1}\right]$ y $\left[\overline{\tau}_{1},\tau_{2}\right]$ son independientes, la FGP conjunta para el n\'umero de usuarios en el sistema al tiempo $t=\tau_{2}$, la FGP conjunta para el n\'umero de usuarios en el sistema est\'a dada por

\begin{eqnarray*}
F_{1}\left(\mathbf{z}\right)&=&R_{2}\left(\prod_{i=1}^{2}P\left(z_{i}\right)\right)F_{2}\left(z_{1},\theta_{2}\left(P_{1}\left(z_{1}\right)\right)\right),\\
F_{2}\left(\mathbf{z}\right)&=&R_{1}\left(\prod_{i=1}^{2}P\left(z_{i}\right)\right)F_{1}\left(\theta_{1}\left(P_{2}\left(z_{2}\right)\right),z_{2}\right).
\end{eqnarray*}

Entonces debemos de determinar las siguientes expresiones:

\begin{eqnarray*}
\begin{array}{cc}
f_{1}\left(1\right)=\frac{\partial F_{1}\left(\mathbf{z}\right)}{\partial z_{1}}|_{\mathbf{z}=1}, & f_{1}\left(2\right)=\frac{\partial F_{1}\left(\mathbf{z}\right)}{\partial z_{2}}|_{\mathbf{z}=1},\\
f_{2}\left(1\right)=\frac{\partial F_{2}\left(\mathbf{z}\right)}{\partial z_{1}}|_{\mathbf{z}=1}, & f_{2}\left(2\right)=\frac{\partial F_{2}\left(\mathbf{z}\right)}{\partial z_{2}}|_{\mathbf{z}=1},\\
\end{array}
\end{eqnarray*}
calculando las derivadas parciales 
\begin{eqnarray*}
\begin{array}{ll}
\frac{\partial R_{1}\left(\mathbf{z}\right)}{\partial
z_{1}}|_{\mathbf{z}=1}=R_{1}^{(1)}\left(1\right)P_{1}^{(1)}\left(1\right),&
\frac{\partial R_{1}\left(\mathbf{z}\right)}{\partial
z_{2}}|_{\mathbf{z}=1}=R_{1}^{(1)}\left(1\right)P_{2}^{(1)}\left(1\right),\\
\frac{\partial R_{2}\left(\mathbf{z}\right)}{\partial
z_{1}}|_{\mathbf{z}=1}=R_{2}^{(1)}\left(1\right)P_{1}^{(1)}\left(1\right),&
\frac{\partial R_{2}\left(\mathbf{z}\right)}{\partial
z_{2}}|_{\mathbf{z}=1}=R_{2}^{(1)}\left(1\right)P_{2}^{(1)}\left(1\right).\\
\end{array}
\end{eqnarray*}

igualando a cero

\begin{eqnarray*}
\begin{array}{ll}
\frac{\partial}{\partial
z_{1}}F_{1}\left(\theta_{1}\left(P_{2}\left(z_{2}\right)\right),z_{2}\right)=0,&
\frac{\partial}{\partial
z_{2}}F_{1}\left(\theta_{1}\left(P_{2}\left(z_{2}\right)\right),z_{2}\right)=\frac{\partial
F_{1}}{\partial z_{2}}+\frac{\partial F_{1}}{\partial
z_{1}}\theta_{1}^{(1)}P_{2}^{(1)}\left(1\right),\\
\frac{\partial}{\partial
z_{1}}F_{2}\left(z_{1},\theta_{2}\left(P_{1}\left(z_{1}\right)\right)\right)=\frac{\partial
F_{2}}{\partial z_{1}}+\frac{\partial F_{2}}{\partial
z_{2}}\theta_{2}^{(1)}P_{1}^{(1)}\left(1\right),&
\frac{\partial}{\partial
z_{2}}F_{2}\left(z_{1},\theta_{2}\left(P_{1}\left(z_{1}\right)\right)\right)=0.
\end{array}
\end{eqnarray*}

Por lo tanto de las dos secciones anteriores se tiene que:

\begin{eqnarray*}
\frac{\partial F_{1}}{\partial z_{1}}&=&\frac{\partial
R_{2}}{\partial z_{1}}|_{\mathbf{z}=1}+\frac{\partial F_{2}}{\partial z_{1}}|_{\mathbf{z}=1}=R_{2}^{(1)}\left(1\right)P_{1}^{(1)}\left(1\right)+f_{2}\left(1\right)+f_{2}\left(2\right)\theta_{2}^{(1)}\left(1\right)P_{1}^{(1)}\left(1\right),\\
\frac{\partial F_{1}}{\partial z_{2}}&=&\frac{\partial
R_{2}}{\partial z_{2}}|_{\mathbf{z}=1}+\frac{\partial F_{2}}{\partial z_{2}}|_{\mathbf{z}=1}=R_{2}^{(1)}\left(1\right)P_{2}^{(1)}\left(1\right),\\
\frac{\partial F_{2}}{\partial z_{1}}&=&\frac{\partial
R_{1}}{\partial z_{1}}|_{\mathbf{z}=1}+\frac{\partial F_{1}}{\partial z_{1}}|_{\mathbf{z}=1}=R_{1}^{(1)}\left(1\right)P_{1}^{(1)}\left(1\right),\\
\frac{\partial F_{2}}{\partial z_{2}}&=&\frac{\partial
R_{1}}{\partial z_{2}}|_{\mathbf{z}=1}+\frac{\partial F_{1}}{\partial z_{2}}|_{\mathbf{z}=1}
=R_{1}^{(1)}\left(1\right)P_{2}^{(1)}\left(1\right)+f_{1}\left(1\right)\theta_{1}^{(1)}\left(1\right)P_{2}^{(1)}\left(1\right).
\end{eqnarray*}

El cual se puede escribir en forma equivalente:
\begin{eqnarray*}
\begin{array}{ll}
f_{1}\left(1\right)=r_{2}\mu_{1}+f_{2}\left(1\right)+f_{2}\left(2\right)\frac{\mu_{1}}{1-\mu_{2}},&
f_{1}\left(2\right)=r_{2}\mu_{2},\\
f_{2}\left(1\right)=r_{1}\mu_{1},&
f_{2}\left(2\right)=r_{1}\mu_{2}+f_{1}\left(2\right)+f_{1}\left(1\right)\frac{\mu_{2}}{1-\mu_{1}}.
\end{array}
\end{eqnarray*}

De donde:
\begin{eqnarray*}
\begin{array}{ll}
f_{1}\left(1\right)=\mu_{1}\left[r_{2}+\frac{f_{2}\left(2\right)}{1-\mu_{2}}\right]+f_{2}\left(1\right)&
f_{2}\left(2\right)=\mu_{2}\left[r_{1}+\frac{f_{1}\left(1\right)}{1-\mu_{1}}\right]+f_{1}\left(2\right).
\end{array}
\end{eqnarray*}

Resolviendo para $f_{1}\left(1\right)$:
\begin{eqnarray*}
f_{1}\left(1\right)&=&r_{2}\mu_{1}+f_{2}\left(1\right)+f_{2}\left(2\right)\frac{\mu_{1}}{1-\mu_{2}}=r_{2}\mu_{1}+r_{1}\mu_{1}+f_{2}\left(2\right)\frac{\mu_{1}}{1-\mu_{2}}\\
&=&\mu_{1}\left(r_{2}+r_{1}\right)+f_{2}\left(2\right)\frac{\mu_{1}}{1-\mu_{2}}=\mu_{1}\left(r+\frac{f_{2}\left(2\right)}{1-\mu_{2}}\right),\\
\end{eqnarray*}

entonces

\begin{eqnarray*}
f_{2}\left(2\right)&=&\mu_{2}\left(r_{1}+\frac{f_{1}\left(1\right)}{1-\mu_{1}}\right)+f_{1}\left(2\right)=\mu_{2}\left(r_{1}+\frac{f_{1}\left(1\right)}{1-\mu_{1}}\right)+r_{2}\mu_{2}=\mu_{2}\left[r_{1}+r_{2}+\frac{f_{1}\left(1\right)}{1-\mu_{1}}\right]=\mu_{2}\left[r+\frac{f_{1}\left(1\right)}{1-\mu_{1}}\right]\\
&=&\mu_{2}r+\mu_{1}\left(r+\frac{f_{2}\left(2\right)}{1-\mu_{2}}\right)\frac{\mu_{2}}{1-\mu_{1}}=\mu_{2}r+\mu_{2}\frac{r\mu_{1}}{1-\mu_{1}}+f_{2}\left(2\right)\frac{\mu_{1}\mu_{2}}{\left(1-\mu_{1}\right)\left(1-\mu_{2}\right)}\\
&=&\mu_{2}\left(r+\frac{r\mu_{1}}{1-\mu_{1}}\right)+f_{2}\left(2\right)\frac{\mu_{1}\mu_{2}}{\left(1-\mu_{1}\right)\left(1-\mu_{2}\right)}=\mu_{2}\left(\frac{r}{1-\mu_{1}}\right)+f_{2}\left(2\right)\frac{\mu_{1}\mu_{2}}{\left(1-\mu_{1}\right)\left(1-\mu_{2}\right)}\\
\end{eqnarray*}
entonces
\begin{eqnarray*}
f_{2}\left(2\right)-f_{2}\left(2\right)\frac{\mu_{1}\mu_{2}}{\left(1-\mu_{1}\right)\left(1-\mu_{2}\right)}&=&\mu_{2}\left(\frac{r}{1-\mu_{1}}\right)\\
f_{2}\left(2\right)\left(1-\frac{\mu_{1}\mu_{2}}{\left(1-\mu_{1}\right)\left(1-\mu_{2}\right)}\right)&=&\mu_{2}\left(\frac{r}{1-\mu_{1}}\right)\\
f_{2}\left(2\right)\left(\frac{1-\mu_{1}-\mu_{2}+\mu_{1}\mu_{2}-\mu_{1}\mu_{2}}{\left(1-\mu_{1}\right)\left(1-\mu_{2}\right)}\right)&=&\mu_{2}\left(\frac{r}{1-\mu_{1}}\right)\\
f_{2}\left(2\right)\left(\frac{1-\mu}{\left(1-\mu_{1}\right)\left(1-\mu_{2}\right)}\right)&=&\mu_{2}\left(\frac{r}{1-\mu_{1}}\right)\\
\end{eqnarray*}
por tanto
\begin{eqnarray*}
f_{2}\left(2\right)&=&\frac{r\frac{\mu_{2}}{1-\mu_{1}}}{\frac{1-\mu}{\left(1-\mu_{1}\right)\left(1-\mu_{2}\right)}}=\frac{r\mu_{2}\left(1-\mu_{1}\right)\left(1-\mu_{2}\right)}{\left(1-\mu_{1}\right)\left(1-\mu\right)}=\frac{\mu_{2}\left(1-\mu_{2}\right)}{1-\mu}r=r\mu_{2}\frac{1-\mu_{2}}{1-\mu}.
\end{eqnarray*}
es decir

\begin{eqnarray}
f_{2}\left(2\right)&=&r\mu_{2}\frac{1-\mu_{2}}{1-\mu}.
\end{eqnarray}

Entonces

\begin{eqnarray*}
f_{1}\left(1\right)&=&\mu_{1}r+f_{2}\left(2\right)\frac{\mu_{1}}{1-\mu_{2}}=\mu_{1}r+\left(\frac{\mu_{2}\left(1-\mu_{2}\right)}{1-\mu}r\right)\frac{\mu_{1}}{1-\mu_{2}}=\mu_{1}r+\mu_{1}r\left(\frac{\mu_{2}}{1-\mu}\right)=\mu_{1}r\left[1+\frac{\mu_{2}}{1-\mu}\right]\\
&=&r\mu_{1}\frac{1-\mu_{1}}{1-\mu}.
\end{eqnarray*}

es decir
\begin{equation}
f_{1}\left(1\right)=r\mu_{1}\frac{1-\mu_{1}}{1-\mu}.
\end{equation}

\section*{Conclusiones}
Los sistemas de visitas son procesos estoc\'asticos en los cuales, mediante el uso de la Funci\'on Generadora de Probabilidades \cite{Adan}, es posible determinar las principales m\'etricas de desempe\~no tanto del sistema en su conjunto como de cada una de las colas que lo conforman. Aunque la revisi\'on de los resultados aqu\'i presentada no pretende ser exhaustiva ni alcanzar una profundidad te\'orica considerable, ofrece una visi\'on general lo suficientemente adecuada para comprender los aspectos t\'ecnicos necesarios y derivar las expresiones mostradas. Para profundizar en estos temas, se recomienda consultar las referencias \cite{BosBoon}, \cite{BoonMeiWinands}, \cite{Boxma}, \cite{LevySidi}, \cite{Roubos}, y \cite{Semenova}.

A partir de los conceptos revisados, quienes deseen profundizar en el estudio de estos procesos pueden explorar los modelos de flujo \cite{Down}, \cite{KaspiMandelbaum}, los cuales proporcionan un enfoque te\'orico m\'as riguroso. Este an\'alisis m\'as detallado \cite{Asmussen}, \cite{Daley68}, \cite{KaspiMandelbaum} permite abordar condiciones de estabilidad y formular propiedades como la ley de los grandes n\'umeros, proporcionando as\'i una perspectiva m\'as completa sobre los sistemas de visitas y sus aplicaciones.

Finalmente, extiendo mi m\'as profundo agradecimiento al Dr. Ra\'ul Montes de Oca Machorro y a la Dra. Patricia Saavedra Barrera por su invaluable supervisi\'on y seguimiento durante la realizaci\'on de este estudio, sin los cuales este trabajo no habr\'ia sido posible.

\section*{Ap\'endice A: El problema de la ruina del jugador}

Supongamos que se tiene un jugador que cuenta con un capital inicial de $\tilde{L}_{0}\geq0$ unidades, esta persona realiza una serie de dos juegos simult\'aneos e independientes de manera sucesiva, dichos eventos son independientes e id\'enticos entre s\'i para cada realizaci\'on.\\

Para $n\geq0$ fijo, la ganancia en el $n$-\'esimo juego es $\tilde{X}_{n}=X_{n}+Y_{n}$ unidades de las cuales se resta una cuota de 1 unidad por cada juego simult\'aneo, es decir, se restan dos unidades por cada juego realizado. En t\'erminos de la teor\'ia de colas puede pensarse como el n\'umero de usuarios que llegan a una cola v\'ia dos procesos de arribo distintos e independientes entre s\'i. Su Funci\'on Generadora de Probabilidades (FGP) est\'a dada por $F\left(z\right)=\esp\left[z^{\tilde{L}_{0}}\right]$ para $z\in\mathbb{C}$, adem\'as
\begin{eqnarray}
\tilde{P}\left(z\right)=\esp\left[z^{\tilde{X}_{n}}\right]=\esp\left[z^{X_{n}+Y_{n}}\right]=\esp\left[z^{X_{n}}z^{Y_{n}}\right]=\esp\left[z^{X_{n}}\right]\esp\left[z^{Y_{n}}\right]=P\left(z\right)\check{P}\left(z\right),
\end{eqnarray}
con $\tilde{\mu}=\esp\left[\tilde{X}_{n}\right]=\tilde{P}\left[z\right]<1$. Sea $\tilde{L}_{n}$ el capital remanente despu\'es del $n$-\'esimo
juego. Entonces

\begin{eqnarray}
\tilde{L}_{n}=\tilde{L}_{0}+\tilde{X}_{1}+\tilde{X}_{2}+\cdots+\tilde{X}_{n}-2n.
\end{eqnarray}

La ruina del jugador ocurre despu\'es del $n$-\'esimo juego, es decir, la cola se vac\'ia despu\'es del $n$-\'esimo juego, entonces sea $T$ definida como $T=min\left\{\tilde{L}_{n}=0\right\}$. Si $\tilde{L}_{0}=0$, entonces claramente $T=0$. En este sentido $T$ puede interpretarse como la longitud del periodo de tiempo que el servidor ocupa para dar servicio en la cola, comenzando con $\tilde{L}_{0}$ grupos de usuarios presentes en la cola, quienes arribaron conforme a un proceso dado por $\tilde{P}\left(z\right)$.\\

Sea $g_{n,k}$ la probabilidad del evento de que el jugador no caiga en ruina antes del $n$-\'esimo juego, y que adem\'as tenga un capital de $k$ unidades antes del $n$-\'esimo juego, es decir, dada $n\in\left\{1,2,\ldots\right\}$ y $k\in\left\{0,1,2,\ldots\right\}$
\begin{eqnarray}
g_{n,k}:=P\left\{\tilde{L}_{j}>0, j=1,\ldots,n,
\tilde{L}_{n}=k\right\},
\end{eqnarray}
la cual adem\'as se puede escribir como:
\begin{eqnarray*}
g_{n,k}&=&P\left\{\tilde{L}_{j}>0, j=1,\ldots,n,
\tilde{L}_{n}=k\right\}=\sum_{j=1}^{k+1}g_{n-1,j}P\left\{\tilde{X}_{n}=k-j+1\right\}\\
&=&\sum_{j=1}^{k+1}g_{n-1,j}P\left\{X_{n}+Y_{n}=k-j+1\right\}=\sum_{j=1}^{k+1}\sum_{l=1}^{j}g_{n-1,j}P\left\{X_{n}+Y_{n}=k-j+1,Y_{n}=l\right\}\\
&=&\sum_{j=1}^{k+1}\sum_{l=1}^{j}g_{n-1,j}P\left\{X_{n}+Y_{n}=k-j+1|Y_{n}=l\right\}P\left\{Y_{n}=l\right\}\\
&=&\sum_{j=1}^{k+1}\sum_{l=1}^{j}g_{n-1,j}P\left\{X_{n}=k-j-l+1\right\}P\left\{Y_{n}=l\right\},
\end{eqnarray*}

es decir
\begin{eqnarray}\label{Eq.Gnk.2S}
g_{n,k}=\sum_{j=1}^{k+1}\sum_{l=1}^{j}g_{n-1,j}P\left\{X_{n}=k-j-l+1\right\}P\left\{Y_{n}=l\right\}.
\end{eqnarray}
Adem\'as
\begin{equation}\label{Eq.L02S}
g_{0,k}=P\left\{\tilde{L}_{0}=k\right\}.
\end{equation}
Se definen las siguientes FGP:
\begin{equation}\label{Eq.3.16.a.2S}
G_{n}\left(z\right)=\sum_{k=0}^{\infty}g_{n,k}z^{k},\textrm{ para
}n=0,1,\ldots,
\end{equation}
y 
\begin{equation}\label{Eq.3.16.b.2S}
G\left(z,w\right)=\sum_{n=0}^{\infty}G_{n}\left(z\right)w^{n}, z,w\in\mathbb{C}.
\end{equation}
En particular para $k=0$,
\begin{eqnarray*}
g_{n,0}=G_{n}\left(0\right)=P\left\{\tilde{L}_{j}>0,\textrm{ para
}j<n,\textrm{ y }\tilde{L}_{n}=0\right\}=P\left\{T=n\right\},
\end{eqnarray*}
adem\'as
\begin{eqnarray*}
G\left(0,w\right)=\sum_{n=0}^{\infty}G_{n}\left(0\right)w^{n}=\sum_{n=0}^{\infty}P\left\{T=n\right\}w^{n}
=\esp\left[w^{T}\right].
\end{eqnarray*}
la cu\'al resulta ser la FGP del tiempo de ruina $T$.

\begin{Prop}\label{Prop.1.1.2S}
Sean $z,w\in\mathbb{C}$ y sea $n\geq0$ fijo. Para $G_{n}\left(z\right)$ y $G\left(z,w\right)$ definidas como en (\ref{Eq.3.16.a.2S}) y (\ref{Eq.3.16.b.2S}) respectivamente, se tiene que
\begin{equation}\label{Eq.Pag.45}
G_{n}\left(z\right)=\frac{1}{z}\left[G_{n-1}\left(z\right)-G_{n-1}\left(0\right)\right]\tilde{P}\left(z\right).
\end{equation}

Adem\'as

\begin{equation}\label{Eq.Pag.46}
G\left(z,w\right)=\frac{zF\left(z\right)-wP\left(z\right)G\left(0,w\right)}{z-wR\left(z\right)},
\end{equation}

con un \'unico polo en el c\'irculo unitario, adem\'as, el polo es
de la forma $z=\theta\left(w\right)$ y satisface que

\begin{enumerate}
\item[i)]$\tilde{\theta}\left(1\right)=1$,

\item[ii)] $\tilde{\theta}^{(1)}\left(1\right)=\frac{1}{1-\tilde{\mu}}$,

\item[iii)]
$\tilde{\theta}^{(2)}\left(1\right)=\frac{\tilde{\mu}}{\left(1-\tilde{\mu}\right)^{2}}+\frac{\tilde{\sigma}}{\left(1-\tilde{\mu}\right)^{3}}$.
\end{enumerate}

Finalmente, adem\'as se cumple que
\begin{equation}
\esp\left[w^{T}\right]=G\left(0,w\right)=F\left[\tilde{\theta}\left(w\right)\right].
\end{equation}

\textbf{Demostraci\'on:}\\
Multiplicando las ecuaciones (\ref{Eq.Gnk.2S}) y (\ref{Eq.L02S}) por el t\'ermino $z^{k}$:

\begin{eqnarray*}
\begin{array}{ll}
g_{n,k}z^{k}=\sum_{j=1}^{k+1}\sum_{l=1}^{j}g_{n-1,j}P\left\{X_{n}=k-j-l+1\right\}P\left\{Y_{n}=l\right\}z^{k},&\textrm{ y  }g_{0,k}z^{k}=P\left\{\tilde{L}_{0}=k\right\}z^{k},

\end{array}
\end{eqnarray*}

ahora sumamos sobre $k$
\begin{eqnarray*}
\sum_{k=0}^{\infty}g_{n,k}z^{k}&=&\sum_{k=0}^{\infty}\sum_{j=1}^{k+1}\sum_{l=1}^{j}g_{n-1,j}P\left\{X_{n}=k-j-l+1\right\}P\left\{Y_{n}=l\right\}z^{k}\\
&=&\sum_{k=0}^{\infty}z^{k}\sum_{j=1}^{k+1}\sum_{l=1}^{j}g_{n-1,j}P\left\{X_{n}=k-\left(j+l-1\right)\right\}P\left\{Y_{n}=l\right\}\\
&=&\sum_{k=0}^{\infty}z^{k+\left(j+l-1\right)-\left(j+l-1\right)}\sum_{j=1}^{k+1}\sum_{l=1}^{j}g_{n-1,j}P\left\{X_{n}=k-\left(j+l-1\right)\right\}P\left\{Y_{n}=l\right\}\\
&=&\sum_{k=0}^{\infty}\sum_{j=1}^{k+1}\sum_{l=1}^{j}g_{n-1,j}z^{j-1}P\left\{X_{n}=k-\left(j+l-1\right)\right\}z^{k-\left(j+l-1\right)}P\left\{Y_{n}=l\right\}z^{l}\\
&=&\sum_{j=1}^{\infty}\sum_{l=1}^{j}g_{n-1,j}z^{j-1}\sum_{k=j+l-1}^{\infty}P\left\{X_{n}=k-\left(j+l-1\right)\right\}z^{k-\left(j+l-1\right)}P\left\{Y_{n}=l\right\}z^{l}\\
&=&\sum_{j=1}^{\infty}g_{n-1,j}z^{j-1}\sum_{l=1}^{j}\sum_{k=j+l-1}^{\infty}P\left\{X_{n}=k-\left(j+l-1\right)\right\}z^{k-\left(j+l-1\right)}P\left\{Y_{n}=l\right\}z^{l}\\
&=&\sum_{j=1}^{\infty}g_{n-1,j}z^{j-1}\sum_{k=j+l-1}^{\infty}\sum_{l=1}^{j}P\left\{X_{n}=k-\left(j+l-1\right)\right\}z^{k-\left(j+l-1\right)}P\left\{Y_{n}=l\right\}z^{l}\\
&=&\sum_{j=1}^{\infty}g_{n-1,j}z^{j-1}\sum_{k=j+l-1}^{\infty}\sum_{l=1}^{j}P\left\{X_{n}=k-\left(j+l-1\right)\right\}z^{k-\left(j+l-1\right)}\sum_{l=1}^{j}P\left\{Y_{n}=l\right\}z^{l}\\
&=&\sum_{j=1}^{\infty}g_{n-1,j}z^{j-1}\sum_{l=1}^{\infty}P\left\{Y_{n}=l\right\}z^{l}\sum_{k=j+l-1}^{\infty}\sum_{l=1}^{j}P\left\{X_{n}=k-\left(j+l-1\right)\right\}z^{k-\left(j+l-1\right)}\\
&=&\frac{1}{z}\left[G_{n-1}\left(z\right)-G_{n-1}\left(0\right)\right]\check{P}\left(z\right)\sum_{k=j+l-1}^{\infty}\sum_{l=1}^{j}P\left\{X_{n}=k-\left(j+l-1\right)\right\}z^{k-\left(j+l-1\right)}\\
&=&\frac{1}{z}\left[G_{n-1}\left(z\right)-G_{n-1}\left(0\right)\right]\check{P}\left(z\right)P\left(z\right)=\frac{1}{z}\left[G_{n-1}\left(z\right)-G_{n-1}\left(0\right)\right]\tilde{P}\left(z\right),
\end{eqnarray*}
es decir la ecuaci\'on (\ref{Eq.3.16.a.2S}) se puede reescribir como
\begin{equation}\label{Eq.3.16.a.2Sbis}
G_{n}\left(z\right)=\frac{1}{z}\left[G_{n-1}\left(z\right)-G_{n-1}\left(0\right)\right]\tilde{P}\left(z\right).
\end{equation}

Por otra parte recordemos la ecuaci\'on (\ref{Eq.3.16.a.2S})
\begin{eqnarray*}
G_{n}\left(z\right)&=&\sum_{k=0}^{\infty}g_{n,k}z^{k},\textrm{ entonces }\frac{G_{n}\left(z\right)}{z}=\sum_{k=1}^{\infty}g_{n,k}z^{k-1},
\end{eqnarray*}

por lo tanto utilizando la ecuaci\'on (\ref{Eq.3.16.a.2Sbis}):

\begin{eqnarray*}
G\left(z,w\right)&=&\sum_{n=0}^{\infty}G_{n}\left(z\right)w^{n}=G_{0}\left(z\right)+\sum_{n=1}^{\infty}G_{n}\left(z\right)w^{n}=F\left(z\right)+\sum_{n=0}^{\infty}\left[G_{n}\left(z\right)-G_{n}\left(0\right)\right]w^{n}\frac{\tilde{P}\left(z\right)}{z}\\
&=&F\left(z\right)+\frac{w}{z}\sum_{n=0}^{\infty}\left[G_{n}\left(z\right)-G_{n}\left(0\right)\right]w^{n-1}\tilde{P}\left(z\right),
\end{eqnarray*}
es decir
\begin{eqnarray*}
G\left(z,w\right)&=&F\left(z\right)+\frac{w}{z}\left[G\left(z,w\right)-G\left(0,w\right)\right]\tilde{P}\left(z\right),
\end{eqnarray*}
entonces
\begin{eqnarray*}
G\left(z,w\right)=F\left(z\right)+\frac{w}{z}\left[G\left(z,w\right)-G\left(0,w\right)\right]\tilde{P}\left(z\right)&=&F\left(z\right)+\frac{w}{z}\tilde{P}\left(z\right)G\left(z,w\right)-\frac{w}{z}\tilde{P}\left(z\right)G\left(0,w\right)\\
&\Leftrightarrow&\\
G\left(z,w\right)\left\{1-\frac{w}{z}\tilde{P}\left(z\right)\right\}&=&F\left(z\right)-\frac{w}{z}\tilde{P}\left(z\right)G\left(0,w\right),
\end{eqnarray*}
por lo tanto,
\begin{equation}
G\left(z,w\right)=\frac{zF\left(z\right)-w\tilde{P}\left(z\right)G\left(0,w\right)}{1-w\tilde{P}\left(z\right)}.
\end{equation}
Ahora $G\left(z,w\right)$ es anal\'itica en $|z|=1$. Sean $z,w$ tales que $|z|=1$ y $|w|\leq1$, como $\tilde{P}\left(z\right)$ es FGP,
\begin{eqnarray*}
|z-\left(z-w\tilde{P}\left(z\right)\right)|<|z|\Leftrightarrow|w\tilde{P}\left(z\right)|<|z|
\end{eqnarray*}
es decir, se cumplen las condiciones del Teorema de Rouch\'e y por tanto, $z$ y $z-w\tilde{P}\left(z\right)$ tienen el mismo n\'umero de ceros en $|z|=1$. Sea $z=\tilde{\theta}\left(w\right)$ la soluci\'on \'unica de $z-w\tilde{P}\left(z\right)$, es decir
\begin{equation}\label{Eq.Theta.w}
\tilde{\theta}\left(w\right)-w\tilde{P}\left(\tilde{\theta}\left(w\right)\right)=0,
\end{equation}
con $|\tilde{\theta}\left(w\right)|<1$. Cabe hacer menci\'on que $\tilde{\theta}\left(w\right)$ es la FGP para el tiempo de ruina cuando $\tilde{L}_{0}=1$. Considerando la ecuaci\'on (\ref{Eq.Theta.w})
\begin{eqnarray*}
0&=&\frac{\partial}{\partial w}\tilde{\theta}\left(w\right)|_{w=1}-\frac{\partial}{\partial w}\left\{w\tilde{P}\left(\tilde{\theta}\left(w\right)\right)\right\}|_{w=1}=\tilde{\theta}^{(1)}\left(w\right)|_{w=1}-\frac{\partial}{\partial w}w\left\{\tilde{P}\left(\tilde{\theta}\left(w\right)\right)\right\}|_{w=1}\\
&-&w\frac{\partial}{\partial w}\tilde{P}\left(\tilde{\theta}\left(w\right)\right)|_{w=1}=\tilde{\theta}^{(1)}\left(1\right)-\tilde{P}\left(\tilde{\theta}\left(1\right)\right)-w\left\{\frac{\partial \tilde{P}\left(\tilde{\theta}\left(w\right)\right)}{\partial \tilde{\theta}\left(w\right)}\cdot\frac{\partial\tilde{\theta}\left(w\right)}{\partial w}|_{w=1}\right\}\\
&=&\tilde{\theta}^{(1)}\left(1\right)-\tilde{P}\left(\tilde{\theta}\left(1\right)
\right)-\tilde{P}^{(1)}\left(\tilde{\theta}\left(1\right)\right)\cdot\tilde{\theta}^{(1)}\left(1\right),
\end{eqnarray*}
luego
$$\tilde{P}\left(\tilde{\theta}\left(1\right)\right)=\tilde{\theta}^{(1)}\left(1\right)-\tilde{P}^{(1)}\left(\tilde{\theta}\left(1\right)\right)\cdot\tilde{\theta}^{(1)}\left(1\right)=\tilde{\theta}^{(1)}\left(1\right)\left(1-\tilde{P}^{(1)}\left(\tilde{\theta}\left(1\right)\right)\right),$$
por tanto $$\tilde{\theta}^{(1)}\left(1\right)=\frac{\tilde{P}\left(\tilde{\theta}\left(1\right)\right)}{\left(1-\tilde{P}^{(1)}\left(\tilde{\theta}\left(1\right)\right)\right)}=\frac{1}{1-\tilde{\mu}}.$$
Ahora determinemos el segundo momento de $\tilde{\theta}\left(w\right)$,
nuevamente consideremos la ecuaci\'on (\ref{Eq.Theta.w}):
\begin{eqnarray*}
0&=&\tilde{\theta}\left(w\right)-w\tilde{P}\left(\tilde{\theta}\left(w\right)\right)\Rightarrow 0=\frac{\partial}{\partial w}\left\{\tilde{\theta}\left(w\right)-w\tilde{P}\left(\tilde{\theta}\left(w\right)\right)\right\}\Rightarrow 0=\frac{\partial}{\partial w}\left\{\frac{\partial}{\partial w}\left\{\tilde{\theta}\left(w\right)-w\tilde{P}\left(\tilde{\theta}\left(w\right)\right)\right\}\right\},
\end{eqnarray*}
luego se tiene
\begin{eqnarray*}
&&\frac{\partial}{\partial w}\left\{\frac{\partial}{\partial w}\tilde{\theta}\left(w\right)-\frac{\partial}{\partial w}\left[w\tilde{P}\left(\tilde{\theta}\left(w\right)\right)\right]\right\}
=\frac{\partial}{\partial w}\left\{\frac{\partial}{\partial w}\tilde{\theta}\left(w\right)-\frac{\partial}{\partial w}\left[w\tilde{P}\left(\tilde{\theta}\left(w\right)\right)\right]\right\}\\
&=&\frac{\partial}{\partial w}\left\{\frac{\partial \tilde{\theta}\left(w\right)}{\partial w}-\left[\tilde{P}\left(\tilde{\theta}\left(w\right)\right)+w\frac{\partial}{\partial w}P\left(\tilde{\theta}\left(w\right)\right)\right]\right\}\\
&=&\frac{\partial}{\partial w}\left\{\frac{\partial \tilde{\theta}\left(w\right)}{\partial w}-\left(\tilde{P}\left(\tilde{\theta}\left(w\right)\right)+w\frac{\partial \tilde{P}\left(\tilde{\theta}\left(w\right)\right)}{\partial w}\frac{\partial \tilde{\theta}\left(w\right)}{\partial w}\right]\right\}\\
&=&\frac{\partial}{\partial w}\left\{\tilde{\theta}^{(1)}\left(w\right)-\tilde{P}\left(\tilde{\theta}\left(w\right)\right)-w\tilde{P}^{(1)}\left(\tilde{\theta}\left(w\right)\right)\tilde{\theta}^{(1)}\left(w\right)\right\}\\
&=&\frac{\partial}{\partial w}\tilde{\theta}^{(1)}\left(w\right)-\frac{\partial}{\partial w}\tilde{P}\left(\tilde{\theta}\left(w\right)\right)-\frac{\partial}{\partial w}\left[w\tilde{P}^{(1)}\left(\tilde{\theta}\left(w\right)\right)\tilde{\theta}^{(1)}\left(w\right)\right]\\
&=&\frac{\partial}{\partial w}\tilde{\theta}^{(1)}\left(w\right)-\frac{\partial\tilde{P}\left(\tilde{\theta}\left(w\right)\right)}{\partial\tilde{\theta}\left(w\right)}\frac{\partial \tilde{\theta}\left(w\right)}{\partial w}-\tilde{P}^{(1)}\left(\tilde{\theta}\left(w\right)\right)\tilde{\theta}^{(1)}\left(w\right)-w\frac{\partial\tilde{P}^{(1)}\left(\tilde{\theta}\left(w\right)\right)}{\partial w}\tilde{\theta}^{(1)}\left(w\right)\\
&-&w\tilde{P}^{(1)}\left(\tilde{\theta}\left(w\right)\right)\frac{\partial \tilde{\theta}^{(1)}\left(w\right)}{\partial w}\\
&=&\tilde{\theta}^{(2)}\left(w\right)-\tilde{P}^{(1)}\left(\tilde{\theta}\left(w\right)\right)\tilde{\theta}^{(1)}\left(w\right)-\tilde{P}^{(1)}\left(\tilde{\theta}\left(w\right)\right)\tilde{\theta}^{(1)}\left(w\right)-w\tilde{P}^{(2)}\left(\tilde{\theta}\left(w\right)\right)\left(\tilde{\theta}^{(1)}\left(w\right)\right)^{2}\\
&-&w\tilde{P}^{(1)}\left(\tilde{\theta}\left(w\right)\right)\tilde{\theta}^{(2)}\left(w\right)\\
&=&\tilde{\theta}^{(2)}\left(w\right)-2\tilde{P}^{(1)}\left(\tilde{\theta}\left(w\right)\right)\tilde{\theta}^{(1)}\left(w\right)-w\tilde{P}^{(2)}\left(\tilde{\theta}\left(w\right)\right)\left(\tilde{\theta}^{(1)}\left(w\right)\right)^{2}-w\tilde{P}^{(1)}\left(\tilde{\theta}\left(w\right)\right)\tilde{\theta}^{(2)}\left(w\right)\\
&=&\tilde{\theta}^{(2)}\left(w\right)\left[1-w\tilde{P}^{(1)}\left(\tilde{\theta}\left(w\right)\right)\right]-
\tilde{\theta}^{(1)}\left(w\right)\left[w\tilde{\theta}^{(1)}\left(w\right)\tilde{P}^{(2)}\left(\tilde{\theta}\left(w\right)\right)+2\tilde{P}^{(1)}\left(\tilde{\theta}\left(w\right)\right)\right],
\end{eqnarray*}
luego
\begin{eqnarray*}
\tilde{\theta}^{(2)}\left(w\right)&&\left[1-w\tilde{P}^{(1)}\left(\tilde{\theta}\left(w\right)\right)\right]-\tilde{\theta}^{(1)}\left(w\right)\left[w\tilde{\theta}^{(1)}\left(w\right)\tilde{P}^{(2)}\left(\tilde{\theta}\left(w\right)\right)+2\tilde{P}^{(1)}\left(\tilde{\theta}\left(w\right)\right)\right]=0\\
\tilde{\theta}^{(2)}\left(w\right)&=&\frac{\tilde{\theta}^{(1)}\left(w\right)\left[w\tilde{\theta}^{(1)}\left(w\right)\tilde{P}^{(2)}\left(\tilde{\theta}\left(w\right)\right)+2P^{(1)}\left(\tilde{\theta}\left(w\right)\right)\right]}{1-w\tilde{P}^{(1)}\left(\tilde{\theta}\left(w\right)\right)}\\
&=&\frac{\tilde{\theta}^{(1)}\left(w\right)w\tilde{\theta}^{(1)}\left(w\right)\tilde{P}^{(2)}\left(\tilde{\theta}\left(w\right)\right)}{1-w\tilde{P}^{(1)}\left(\tilde{\theta}\left(w\right)\right)}+\frac{2\tilde{\theta}^{(1)}\left(w\right)\tilde{P}^{(1)}\left(\tilde{\theta}\left(w\right)\right)}{1-w\tilde{P}^{(1)}\left(\tilde{\theta}\left(w\right)\right)}.
\end{eqnarray*}
Si evaluamos la expresi\'on anterior en $w=1$:
\begin{eqnarray*}
\tilde{\theta}^{(2)}\left(1\right)&=&\frac{\left(\tilde{\theta}^{(1)}\left(1\right)\right)^{2}\tilde{P}^{(2)}\left(\tilde{\theta}\left(1\right)\right)}{1-\tilde{P}^{(1)}\left(\tilde{\theta}\left(1\right)\right)}+\frac{2\tilde{\theta}^{(1)}\left(1\right)\tilde{P}^{(1)}\left(\tilde{\theta}\left(1\right)\right)}{1-\tilde{P}^{(1)}\left(\tilde{\theta}\left(1\right)\right)}=\frac{\left(\tilde{\theta}^{(1)}\left(1\right)\right)^{2}\tilde{P}^{(2)}\left(1\right)}{1-\tilde{P}^{(1)}\left(1\right)}+\frac{2\tilde{\theta}^{(1)}\left(1\right)\tilde{P}^{(1)}\left(1\right)}{1-\tilde{P}^{(1)}\left(1\right)}\\
&=&\frac{\left(\frac{1}{1-\tilde{\mu}}\right)^{2}\tilde{P}^{(2)}\left(1\right)}{1-\tilde{\mu}}+\frac{2\left(\frac{1}{1-\tilde{\mu}}\right)\tilde{\mu}}{1-\tilde{\mu}}=\frac{\tilde{P}^{(2)}\left(1\right)}{\left(1-\tilde{\mu}\right)^{3}}+\frac{2\tilde{\mu}}{\left(1-\tilde{\mu}\right)^{2}}=\frac{\sigma^{2}-\tilde{\mu}+\tilde{\mu}^{2}}{\left(1-\tilde{\mu}\right)^{3}}+\frac{2\tilde{\mu}}{\left(1-\tilde{\mu}\right)^{2}}\\
&=&\frac{\sigma^{2}-\tilde{\mu}+\tilde{\mu}^{2}+2\tilde{\mu}\left(1-\tilde{\mu}\right)}{\left(1-\tilde{\mu}\right)^{3}},
\end{eqnarray*}
es decir,
\begin{eqnarray*}
\tilde{\theta}^{(2)}\left(1\right)&=&\frac{\sigma^{2}}{\left(1-\tilde{\mu}\right)^{3}}+\frac{\tilde{\mu}}{\left(1-\tilde{\mu}\right)^{2}}.
\end{eqnarray*}
\end{Prop}

\begin{Cor}
El tiempo de ruina del jugador tiene primer y segundo momento dados por
\begin{eqnarray}
\esp\left[T\right]&=&\frac{\esp\left[\tilde{L}_{0}\right]}{1-\tilde{\mu}}\\
Var\left[T\right]&=&\frac{Var\left[\tilde{L}_{0}\right]}{\left(1-\tilde{\mu}\right)^{2}}+\frac{\sigma^{2}\esp\left[\tilde{L}_{0}\right]}{\left(1-\tilde{\mu}\right)^{3}}.
\end{eqnarray}
\end{Cor}

\section*{Ap\'endice B: Extensi\'on a sistema de visitas c\'iclicas}
Hagamos una extensi\'on para el caso en que se tienen varias colas:

\begin{Prop}
Supongamos

\begin{equation}\label{Eq.1}
f_{i}\left(i\right)-f_{j}\left(i\right)=\mu_{i}\left[\sum_{k=j}^{i-1}r_{k}+\sum_{k=j}^{i-1}\frac{f_{k}\left(k\right)}{1-\mu_{k}}\right]
\end{equation}

\begin{equation}\label{Eq.2}
f_{i+1}\left(i\right)=r_{i}\mu_{i},
\end{equation}

Demostrar que

\begin{eqnarray*}
f_{i}\left(i\right)&=&\mu_{i}\left[\sum_{k=1}^{N}r_{k}+\sum_{k=1,k\neq i}^{N}\frac{f_{k}\left(k\right)}{1-\mu_{k}}\right].
\end{eqnarray*}

En la Ecuaci\'on (\ref{Eq.2}) hagamos $j=i+1$, entonces se tiene $f_{j}=r_{i}\mu_{i}$, lo mismo para (\ref{Eq.1})

\begin{eqnarray*}
f_{i}\left(i\right)&=&r_{i}\mu_{i}+\mu_{i}\left[\sum_{k=j}^{i-1}r_{k}+\sum_{k=j}^{i-1}\frac{f_{k}\left(k\right)}{1-\mu_{k}}\right]=\mu_{i}\left[\sum_{k=j}^{i}r_{k}+\sum_{k=j}^{i-1}\frac{f_{k}\left(k\right)}{1-\mu_{k}}\right]
\end{eqnarray*}

entonces, tomando sobre todo valor de $1,\ldots,N$, tanto para antes de $i$ como para despu\'es de $i$, entonces

\begin{eqnarray*}
f_{i}\left(i\right)&=&\mu_{i}\left[\sum_{k=1}^{N}r_{k}+\sum_{k=1,k\neq i}^{N}\frac{f_{k}\left(k\right)}{1-\mu_{k}}\right].
\end{eqnarray*}
\end{Prop}

Ahora, supongamos nuevamente la ecuaci\'on (\ref{Eq.1})

\begin{eqnarray*}
f_{i}\left(i\right)-f_{j}\left(i\right)&=&\mu_{i}\left[\sum_{k=j}^{i-1}r_{k}+\sum_{k=j}^{i-1}\frac{f_{k}\left(k\right)}{1-\mu_{k}}\right]\Leftrightarrow
f_{j}\left(j\right)-f_{i}\left(j\right)=\mu_{j}\left[\sum_{k=i}^{j-1}r_{k}+\sum_{k=i}^{j-1}\frac{f_{k}\left(k\right)}{1-\mu_{k}}\right]\\
f_{i}\left(j\right)&=&f_{j}\left(j\right)-\mu_{j}\left[\sum_{k=i}^{j-1}r_{k}+\sum_{k=i}^{j-1}\frac{f_{k}\left(k\right)}{1-\mu_{k}}\right]=\mu_{j}\left(1-\mu_{j}\right)\frac{r}{1-\mu}-\mu_{j}\left[\sum_{k=i}^{j-1}r_{k}+\sum_{k=i}^{j-1}\frac{f_{k}\left(k\right)}{1-\mu_{k}}\right]\\
&=&\mu_{j}\left[\left(1-\mu_{j}\right)\frac{r}{1-\mu}-\sum_{k=i}^{j-1}r_{k}-\sum_{k=i}^{j-1}\frac{f_{k}\left(k\right)}{1-\mu_{k}}\right]=\mu_{j}\left[\left(1-\mu_{j}\right)\frac{r}{1-\mu}-\sum_{k=i}^{j-1}r_{k}-\frac{r}{1-\mu}\sum_{k=i}^{j-1}\mu_{k}\right]\\
&=&\mu_{j}\left[\frac{r}{1-\mu}\left(1-\mu_{j}-\sum_{k=i}^{j-1}\mu_{k}\right)-\sum_{k=i}^{j-1}r_{k}\right]=\mu_{j}\left[\frac{r}{1-\mu}\left(1-\sum_{k=i}^{j}\mu_{k}\right)-\sum_{k=i}^{j-1}r_{k}\right].
\end{eqnarray*}

Ahora,

\begin{eqnarray*}
1-\sum_{k=i}^{j}\mu_{k}=1-\sum_{k=1}^{N}\mu_{k}+\sum_{k=j+1}^{i-1}\mu_{k}\Leftrightarrow
\sum_{k=i}^{j}\mu_{k}=\sum_{k=1}^{N}\mu_{k}-\sum_{k=j+1}^{i-1}\mu_{k}\Leftrightarrow\sum_{k=1}^{N}\mu_{k}=\sum_{k=i}^{j}\mu_{k}+\sum_{k=j+1}^{i-1}\mu_{k}.
\end{eqnarray*}

Por tanto
\begin{eqnarray}
f_{i}\left(j\right)&=&\mu_{j}\left[\frac{r}{1-\mu}\sum_{k=j+1}^{i-1}\mu_{k}+\sum_{k=j}^{i-1}r_{k}\right].
\end{eqnarray}


\end{document}